\documentclass[12pt,leqno]{article}
\usepackage{amssymb,amsfonts,amsmath,amsthm,amscd,mathrsfs}
\def\R{{{{\rm l} \kern -.15em {\rm R}}}}
\def\N{{{{\rm l} \kern -.15em {\rm N}}}}
\def\E{{{{\rm l} \kern -.15em {\rm E}}}}
\def\P{{{{\rm l} \kern -.15em {\rm P}}}}
\def\D{{{{\rm l} \kern -.15em {\rm D}}}}
\def\L{{{{\rm l} \kern -.15em {\rm L}}}}
\def\Z{{{{\rm Z} \kern -.35em {\rm Z}}}}

\def\IR{{\mathbb R}}
\def\IN{{\mathbb N}}

\def\IZ{{\mathbb Z}}

\def\hs{\hspace{-1ex}}
\def\n{\noindent}
\def\dsl{\textstyle\sum\limits}

\def\dis{\displaystyle}

\def\fr{\mbox{\footnotesize $\dis\frac{1}{2}$}}

\def\ov{\overline}
\def\e{\epsilon}
\def\f{\footnotesize}

\def\point{{\mbox{\large $.$}}}
\def\wh{\widehat}
\def\wt{\widetilde}

\def\cA{{\cal A}}

\def\cC{{\cal C}}
\def\cD{{\cal D}}
\def\cE{{\cal E}}
\def\cL{{\cal L}}

\def\cE{{\cal E}}

\def\cK{{\cal K}}

\def\cF{{\cal F}}
\def\cR{{\cal R}}
\def\cS{{\cal S}}
\def\cG{{\cal G}}

\def\cU{{\cal U}}
\def\cV{{\cal V}}
\def\cW{{\cal W}}

\dimendef\dimen=0

\newtheorem{theorem}{Theorem}[section]
\newtheorem{lemma}[theorem]{Lemma}
\newtheorem{corollary}[theorem]{Corollary}
\newtheorem{proposition}[theorem]{Proposition}
\newtheorem{remark}[theorem]{Remark}

\begin{document}

\noindent

~

\bigskip
\begin{center}
{\bf GIANT COMPONENT AND VACANT SET FOR \\RANDOM  WALK ON A DISCRETE TORUS}
\end{center}

%\vspace{1cm}
\begin{center}
Itai Benjamini\footnote[1]{Department of Mathematics, The Weizmann Institute of Science,
POB 26, Rehovot 76100, Israel} and Alain-Sol Sznitman\footnote[2]{Departement Mathematik, ETH Zurich, 
CH-8092 Z\"urich, Switzerland} 
\end{center}

\bigskip
\begin{abstract}
We consider random walk on a discrete torus $E$ of side-length $N$, in sufficiently high dimension $d$. We investigate the percolative properties of the vacant set corresponding to the collection of sites which have not been visited by the walk up to time $uN^d$. We show that when $u$ is chosen small, as $N$ tends to infinity, there is with overwhelming probability a unique connected component in the vacant set which contains segments of length const $\log N$. Moreover, this connected component occupies a non-degenerate fraction of the total number of sites $N^d$ of $E$, and any point of $E$ lies within distance $N^\beta$ of this component, with $\beta$ an arbitrary positive number.
\end{abstract}

 \vfill

\newpage

\thispagestyle{empty}
~

\newpage
\setcounter{page}{1}

\setcounter{section}{-1}
 \section{Introduction}

 \setcounter{equation}{0}

We investigate here random walk on a $d$-dimensional torus of large side-length $N$, and are interested in the set of points that have not been visited by the walk up to times of order $N^d$. This time scale is much shorter than the typical time it takes the walk to cover the discrete torus. Indeed the cover time of the discrete torus is known to be of order $N^d \log N$, when $d \ge 3$, and $N^2 (\log N)^2$, when $d=2$, cf.~\cite{Aldo83}, \cite{AldoFill99}, \cite{BrumHilh91}, \cite{DembPereRoseZeit04}, \cite{DembPereRoseZeit06}, and the references therein, for this and much more. In fact, when $d \ge 3$, and $u$ is an arbitrary positive number, the probability that the walk visits a given point of the discrete torus up to time $uN^d$ remains bounded away from $0$ and $1$ as $N$ tends to infinity. This makes the time scale $N^d$ an appropriate choice to discuss the percolative properties of the vacant set left by the walk. Incidentally these questions are closely related to the analysis of the disconnection time of a discrete cylinder by a random walk, recently investigated in \cite{DembSzni06},\cite{Szni06}. The main object of this work is to discuss the typical presence of a well-defined giant connected component in the vacant set left by the walk by time $uN^d$, for large $N$, when the dimension $d$ is large enough, and $u$ suitably small. We expect a different behavior when $u$ is large, but this work does not present results in this direction. We believe some of our methods and results pertinent to improve bounds of the disconnection time of a discrete cylinder by a random walk recently derived in \cite{DembSzni06}, see also \cite{Szni06}.

\medskip
Before discussing results any further, we describe the model in more details. We consider $d \ge 3$, $N \ge 1$, and denote with $E$ the $d$-dimensional discrete torus of side-length $N$:
\begin{equation}\label{0.1}
E = (\IZ / N\IZ)^d\,.
\end{equation}

\medskip\n
We write $P$, resp. $P_x$, when $x \in E$, for the law on $E^\IN$ endowed with the product $\sigma$-algebra $\cF$, of simple random walk on $E$ started with the uniform distribution, resp. at $x$. We let $X_\point$ stand for the canonical process on $E^\IN$, $X_{[0,t]}$ for the set of sites visited by the walk up to time $[t]$:
\begin{equation}\label{0.2}
\mbox{$X_{[0,t]} = \{z \in E$; for some $0 \le n \le t$, $X_n = z\}$, for $t \ge 0$}\,.
\end{equation}

\medskip\n
Our main focus lies in the percolative properties of the vacant set $E \backslash X_{[0,uN^d]}$  left by the walk up to time $uN^d$, when $N$ is large and $u > 0$ some fixed positive number. We show that when $d \ge 4$ and $u$ suitably small, the vacant set by time $uN^d$, typically contains a profusion of segments of logarithmic size in $N$, for large $N$. More precisely, we define for $K > 0$, $0 < \beta < 1$, $t \ge 0$, the event which specifies that for every point of $E$ there is in each coordinate direction, within $N^\beta$ steps, a segment of length $[K \log N]$ in the vacant set left by the walk at time $t$:
\begin{align}
\cV_{K,\beta,t} = \{ & \mbox{for all $x \in E$, $1 \le j \le d$, for some $0 \le m < N^\beta$,}  \label{0.3}
\\
&X_{[0,t]} \cap \big\{x + \big( m + \big[0, [K \log N]\big]\big) \,e_j\big\} = \phi\big\}\,,\nonumber
\end{align} 

\n
where $(e_i)_{1 \le i \le d}$ stands for the canonical basis of $\IR^d$. We show in Theorem \ref{theo1.2} that for $d,K,\beta$ as above,
\begin{equation}\label{0.4}
\lim\limits_N \;P[\cV_{K,\beta,uN^d}] = 1, \;\mbox{for small $u > 0$}\,.
\end{equation}

\n
We also show in Proposition \ref{prop1.1} that when $d \ge 3$, for $u > 0$, 
\begin{equation}\label{0.5}
e^{-cu} \le \liminf\limits_N \,P[0 \notin X_{[0,uN^d]}] \le \limsup\limits_N \,P[0 \notin X_{[0,uN^d]}] \le e^{-c^\prime u}\,,
\end{equation}

\n
with $c,c^\prime$ suitable positive dimension dependent constants, (more is known, see \cite{AldoFill99}, Chapter 3, Proposition 20 and Chapter 13, Proposition 8). This feature motivates the interest of the time scale $N^d$ in the investigation of the vacant set left by the walk. We sharpen this result by showing in Corollary \ref{cor4.5} that
\begin{equation}\label{0.6}
\lim\limits_N \;P[e^{-cu} \le | \,E \backslash X_{[0,uN^d]}\,| / N^d \le e^{-c^\prime u}] = 1, \;\mbox{for $u > 0$}\,,
\end{equation}

\n
where for $A \subseteq E$, $|A|$ denotes the cardinality of $A$.

\medskip
When the dimension $d$ is suitably large, i.e. $d \ge d_0$, cf.~(\ref{2.39}), we introduce a dimension dependent constant $c_0$, and events $\cG_{\beta,t} \subseteq \cV_{c_0,\beta,t}$, increasing with $\beta \in (0,1)$, such that for any such $\beta$:
\begin{align}
{\rm i)} &\;\; \lim\limits_N \;P[\cG_{\beta,uN^d}] = 1, \;\mbox{for small $u > 0$, and}\nonumber
\\[-1ex]
\label{0.7}
\\[-1ex]
{\rm ii)} &\;\;\; \mbox{on $\cG_{\beta,t}$ there is a unique connected component $O$ in $E \backslash X_{[0,t]}$} \nonumber
\\
&\;\;\, \mbox{which contains segments of length $L_0 = [c_0 \log N]$}\,,\nonumber
\end{align}

\n
(see (\ref{2.51}) for the more general claim). The connected component $O$ is thus well-defined on the nested events $\cG_{\beta,t}$. In view of (\ref{0.3}) and since $\cG_{\beta,t} \subseteq \cV_{c_0, \beta,t}$, the connected component $O$ is ubiquitous on $E$. We refer to it as the {\it giant component}. We also show in Corollary \ref{cor4.6}, that when $u$ is small, $O$ typically has a non-degenerate volume in $E$. More precisely we prove that for $d \ge d_0$, $\beta, \gamma \in (0,1)$,
\begin{equation}\label{0.8}
\lim\limits_N \,P\big[\cG_{\beta,uN^d} \cap \{|O| \ge \gamma N^d\}\big] = 1, \;\mbox{when $u > 0$ is small.}
\end{equation}

\n
However our results do not rule out the possible existence of other components of the vacant set with non-degenerate volume as well, cf.~Remark \ref{rem4.7}. In fact the present work raises many questions. How do percolative properties of the vacant set compare to the picture stemming from Bernoulli bond-percolation? Is there a small $u$ regime with typically one single giant component and all other components of small volume and size, a large $u$ regime with only small connected components and in between a critical regime, (see for instance \cite{HeydHofs05} and references therein)? Simulations performed when $d=3,4,5,6,7,$ seem to support this picture, with a critical threshold located near $u=3$. If such a critical regime can be extracted, do components in the vacant set in this regime inherit some of the invariance properties of Brownian motion viewed as a scaling limit of simple random walk? What are the relevant values of the dimension $d$? It is maybe instructive to also consider these problems on other graphs, such as expanders, (where a small $u$  regime with some giant component  and a large $u$ regime with only small components, can easily be established), random $d$-regular graphs, hypercubes etc., see \cite{AlonBenjStac04} for a study of percolation on such graphs. These are just a few examples of the many questions raised by the present article.

\medskip
We now try to describe some of the ideas and methods involved in the proof of (\ref{0.4}), (\ref{0.7}), (\ref{0.8}).

\medskip
Behind (\ref{0.4}) lies a type of coupon-collector heuristics. We show in Proposition \ref{prop1.1} that up to time $uN^d$ about const $uN^{d-2}$ excursions in and out of two concentric balls centered at the origin with radius some suitable fraction of $N$, take place. At most const $u N^\beta$ of these excursions hit a given segment of length $N^\beta$ starting at the origin. Chopping this segment into $M = N^\beta / [K \log N]$ segments of length $[K \log N]$, and neglecting the possible hits of more than one segment by one such excursion, a coupon-collector heuristics, cf.~Durrett \cite{Durr91}, Chapter 2, Example 6.6, makes it plausible that it takes about $M \log M \sim \frac{\beta}{K} \;N^\beta$ such excursions to hit each of these segments. However when $u$ is chosen small, ${\rm const}\, uN^\beta \ll \frac{\beta}{K} \,N^\beta$, and not all segments can be hit by the walk up to time in $N^d$. The above lines describe the intuition behind the proof of (\ref{0.4}) in Theorem \ref{theo1.2}.

\medskip
The key to the uniqueness statement contained in (\ref{0.7}) is an exponential estimate proved in Theorem \ref{theo2.1}. It shows in particular that when $d \ge 5$, and $\lambda  > 0$ is such that:
\begin{equation}\label{0.9}
e^{2 \lambda} \Big(\mbox{\f $\dis\frac{2}{d}$} + \Big(1 - \mbox{\f $\dis\frac{2}{d}$}\Big) \;q(d-2)\Big) < 1\,,
\end{equation}

\n
where for any integer $\nu \ge 1$,
\begin{equation}\label{0.10}
\mbox{$q(\nu)$ = the return probability to the origin of simple random walk on $\IZ^\nu$}\,,
\end{equation}
then for $N \ge N(d,\lambda)$ and $u \le u(d,\lambda)$,
\begin{equation}\label{0.11}
P[X_{[0,uN^d]} \supseteq A] \le \exp\{-\lambda \,|A|\}\,,
\end{equation}

\n
for any subset $A$ of $E$ contained in the canonical projection $F$ on $E$ of a two-dimensional affine plane generated by two coordinate directions in $\IZ^d$. When $e^\lambda > 7$ can be achieved, (this is  the requirement which specifies $d_0$, cf.~(\ref{2.39})), the exponential bound (\ref{0.11}) combined with a Peierl-type argument yields in Corollary \ref{cor2.5} the key uniqueness statement behind (\ref{0.7}). The claim (\ref{0.7}) is then proved in Corollary \ref{cor2.6}. We also explain in Remark \ref{rem2.4} why a restriction on the class of sets $A$ that appear in (\ref{0.11}) is needed. There is an independent interest to the above exponential bound: a variation of it and of (\ref{0.4}) should lead to a sharpening of the lower bounds on the disconnection time of discrete cylinders $(\IZ/N\IZ)^d \times \IZ$ derived in \cite{DembSzni06}, at least when $d$ is large enough.

\medskip
To prove (\ref{0.8}), we in essence control fluctuations of the proportion of sites in $E$ which at time $uN^d$ are connected by a vacant path in some two-dimensional $F$, as below (\ref{0.11}), to sites at distance $[c_0 \log N]$. Such sites belong to the giant component $O$, when $\cG_{\beta,uN^d}$ occurs, cf.~(\ref{2.54}). This leads us to develop estimates on the covariance of ``local functions'' of the vacant sites left by the walk up to time of order $uN^d$ in the neighborhood of two sufficiently distant point on the torus, see Proposition \ref{prop4.2}. Qualitatively similar issues appear for instance in \cite{DembPereRoseZeit06}. To this end we develop in Theorem \ref{theo3.1} a bound on the total variation norm between the law of a suitable ``limit model'' and the law $Q_{u,w}$ of a recentered excursion of the walk. This excursion runs from the time of the first, up to the last, visit to $C(x_1) \cup C(x_2)$, where $C(x_i)$ are boxes of side-length $2L$ centered at $x_i, i = 1,2$, in $E$, with mutual distance at least $2r + 3$, where $r \ge 10 L$, and the walk is conditioned to start at a point $u$ at distance at least $r$ from $\{x_1,x_2\}$ and exit the $r$-neighborhood of $\{x_1,x_2\}$ at the point $w$, and stop there. The ``limit model'' with law $Q$ corresponds to excursions of the simple random walk on $\IZ^d$ starting with the normalized harmonic measure viewed from infinity of the box $C$ centered at the origin with side-length $2L$, stopped at its last visit of $C$. In Theorem \ref{theo3.1} we show that
\begin{equation}\label{0.12}
\|Q_{u,w} - Q \|_{TV} \le c\,\dis\frac{L^2}{r} \;,
\end{equation}

\n
where $c$ is a dimension dependent constant and $\|\cdot\|_{TV}$ the total variation norm. This estimate is of independent interest and can straightforwardly be extended to the case of finitely many points $x_i$, cf.~Remark \ref{rem3.2}. Our main control on fluctuations of spatial averages on $E$ of local functions is then stated in Theorem \ref{theo4.3}, and enables to show (\ref{0.6}) in Corollary \ref{cor4.5}, and (\ref{0.8}) in Corollary \ref{cor4.6}. In Corollary \ref{cor4.8} we also show that when $d \ge 3$, the largest cube contained in the vacant set at time $uN^d$ typically has size of order $(\log N)^{\frac{1}{d-2}}$, for large $N$. This should be contrasted with the case of Bernoulli bond-percolation on the torus, where for large $N$ the largest cube contained in a cluster typically has much smaller size of order $(\log N)^{\frac{1}{d}}$.  

\medskip
Let us now describe the organization of this article. 

\medskip
In Section 1, we introduce some further notation, and mainly provide the proof of (\ref{0.4}) in Theorem \ref{theo1.2}. Incidentally we show (\ref{0.5}) in Proposition \ref{prop1.1}.

\medskip
In Section 2 we prove a more general version of (\ref{0.11}) in Theorem \ref{theo2.1}, and use it in Corollary \ref{cor2.5}, \ref{cor2.6} to prove (\ref{0.7}), and thereby construct the giant component in the vacant set, which is shown to be typical when $d \ge d_0$ and $u$ small enough.

\medskip
In Section 3, we obtain the total variation estimate (\ref{0.12}) in Theorem \ref{theo3.1}. This comes as a preparation for the control of fluctuations of certain spatial averages of local functions in the next section.

\medskip
In Section 4, we show (\ref{0.8}) in Corollary \ref{cor4.6}, the simpler (\ref{0.6}) in Corollary \ref{cor4.5}, and the controls on the largest cube contained in the vacant set in Corollary \ref{cor4.8}. The variance bounds of Proposition \ref{prop4.2} make strong use of Theorem \ref{theo3.1}. Our general control on fluctuations of averages of local functions appears in Theorem \ref{theo4.3}.

\medskip
Finally throughout the text $c$, or $c^\prime$ denote positive constants which solely depend on $d$, with values that change from place to place. The numbered constants $c_0,c_1,\dots$, are fixed and refer to the value at their first place of appearance in the text. Dependence of constants on additional parameters appears in the notation. For instance $c(K,\beta)$ denotes a positive constant depending on $d,K,\beta$.

\bigskip\bigskip\n
{\bf Acknowledgements:} We wish to thank Laurent Goergen and Roey Izkovsky for their simulations, as well as David Windisch for his comments on a previous version of this work. Alain-Sol Sznitman wants to thank Amir Dembo for helpful discussions.

\section{Ubiquity of vacant segments of logarithmic size}
 \setcounter{equation}{0}
 
 The main object of this section is to show that when $d \ge 4$, for large $N$, up to times that are small multiples of $N^d$ the vacant set left by the walk on the discrete torus $E$ contains with overwhelming probability segments of size of order $\log N$ in the vicinity of each point of $E$, cf.~Theorem \ref{theo1.2} and (\ref{0.4}). We also prove the estimate (\ref{0.5}) on the probability that a point belongs to the vacant set up to time $u N^d$, with $d \ge 3$, cf.~Proposition \ref{prop1.1}. We first need some additional notation.
 
\medskip
We denote with $|\cdot |$ and $|\cdot |_\infty$ the Euclidean and $\ell^\infty$-distances on $\IZ^d$,  or the corresponding distances on $E$. We write $B(x,r)$, for the closed ball relative to $|\cdot|_\infty$, with radius $r \ge 0$, and center $x \in \IZ^d$, or $E$. We denote with $S(x,r)$ the corresponding $|\cdot|_\infty$-sphere with radius $r$ and center $x$. We say that $x,y$ in $\IZ^d$ or $E$ are neighbors, resp.~$\star$-neighbors, if $|x-y| = 1$, resp. $|x-y|_\infty = 1$. The notions of connected or~$\star$-connected subsets of $\IZ^d$ or $E$ are then defined accordingly, and so are the notions of nearest neighbor path or $\star$-nearest neighbor path on $\IZ^d$ or $E$. For $A,B$ subsets of $\IZ^d$ or $E$, we denote with $A+B$ the subset of points of the form  $x+y$, with $x \in A$, $y \in B$. When $U$ is a subset of $\IZ^d$ or $E$, we let $|U|$ stand for the cardinality of $U$ and $\partial U$ for the boundary of $U$:
\begin{equation}\label{1.1}
\partial U = \{x \in U^c; \;\exists y \in U, \;|x - y| = 1\}\,.
\end{equation}

\n
We denote with $\pi_E$ the canonical projection from $\IZ^d$ onto $E$. For $1 \le m \le d$, we write $\cL_m$ for the collection of subsets of $E$ that are projections under $\pi_E$ of affine lattices in $\IZ^d$ generated by $m$ distinct vectors of the canonical basis:
\begin{align}
\cL_m = \Big\{ &\mbox{$F \subseteq E$; for some $I \subseteq \{1,\dots,d\}$  with} \label{1.2}
\\[-1ex]
&| I | = m, \;\;\mbox{and some}\; y \in \IZ^d, \;F = \pi_E \Big(y + \dsl_{i \in I} \IZ \,e_i\Big)\Big\}\,,\nonumber
\end{align}
where as below (\ref{0.3}), $(e_i)_{1 \le i \le d}$ denotes the canonical basis of $\IR^d$.

\medskip
We let $(\theta_n)_{n \ge 0}$ and $(\cF_n)_{n \ge 0}$ stand for the canonical shift on $E^{\IN}$ and the filtration of the canonical process. For $U \subseteq E$, $H_U$ and $T_U$ stand for the entrance time and exit time in or from $U$:
\begin{equation}\label{1.3}
H_U = \inf\{n \ge 0; \;X_n \in U\}, \;\;T_U = \inf\{n \ge 0; \,X_n \notin U\}\,.
\end{equation}
We write $\wt{H}_U$ for the hitting time of $U$:
\begin{equation}\label{1.4}
\wt{H}_U = \inf\{n \ge 1; \;X_n \in U\}\,.
\end{equation}

\n
When $U = \{x\}$, we write as a subscript $x$ in place of $\{x\}$, for simplicity. Given $A \subseteq \wt{A} \subseteq E$, we often consider the successive return times to $A$ and departures from $\wt{A}$:
\begin{equation}\label{1.5}
\begin{array}{l}
R_1  = H_A, \;D_1 = T_{\wt{A}} \circ \theta_{R_1} + R_1, \;\;\mbox{and for $k \ge 1$}
\\[1ex]
R_{k+1}  = H_A \circ \theta_{D_k} + D_k, \;D_{k+1} = D_1 \circ \theta_{D_k} + D_k, \;\;\mbox{so that}
\\[1ex]
0 \le R_1 \le D_1 \le \dots \le R_k \le D_k \le \dots \le \infty
\end{array}
\end{equation}

\medskip\n
and $P$-a.s. the above inequalities are strict except maybe for the first one. We also set $R_0 = 0 = D_0$ by convention. The transition density of the walk on $E$ is denoted by
\begin{equation}\label{1.6}
p_k(x,y) = P_x[X_k = y], \;k \ge 0, \,x,y \in E\,.
\end{equation}

\medskip\n
We write $P_x^{\IZ^\nu}$, or $E^{\IZ^\nu}_x$, for $x \in \IZ^\nu$, $\nu \ge 1$, to refer to the law or expectation for simple random walk on $\IZ^\nu$ starting from $x$. We otherwise keep the same notation as above. We let $g_\nu(\cdot)$ stand for the Green function of simple random walk on $\IZ^\nu$, $\nu \ge 1$, with a pole at the origin:
\begin{equation}\label{1.7}
g_\nu(z) = E_z^{\IZ^\nu} \Big[\dsl_{n \ge 0} \,1_{\{X_n = 0\}}\Big], \;\;\mbox{for $z \in \IZ^\nu$}\,,
\end{equation}

\n
(which of course is identically infinite unless $\nu \ge 3$). As a direct consequence of the geometric number of returns of the walk to the origin, one classically has:
\begin{equation}\label{1.8}
g_\nu(0) = (1 - q(\nu))^{-1}\,,
\end{equation}

\n
where, cf.~(\ref{0.10}), $q(\nu) = P_0^{\IZ^\nu} [\wt{H}_0 < \infty]$, denotes the return probability to the origin.

\medskip
We are now ready to begin and consider for $N \ge 1$
\begin{equation}\label{1.9}
B = \pi_E \Big(\Big[ - \mbox{\f $\dis \frac{N}{8}, \;\frac{N}{8}$}\Big]^d \cap \IZ^d\Big) \subseteq \wt{B} = \pi_E \Big(\Big[ - \mbox{\f $\dis \frac{N}{4}, \;\frac{N}{4}$}\Big]^d \cap \IZ^d\Big) \,,
\end{equation}

\n
as well as, cf.~(\ref{1.5}),
\begin{equation}\label{1.10}
\mbox{$\cR_k$, $\cD_k$, $k \ge 1$, the successive returns to $B$ and departures from $\wt{B}$}\,.
\end{equation}

\n
The following estimates will be useful in the sequel. We also prove the controls (\ref{0.5}) on the probability that a point belongs to the vacant set, see also (2.26) of \cite{BrumHilh91}.

\begin{proposition}\label{prop1.1} $(d \ge 3$)
\begin{align}
&P[\cR_{k^*} \le u N^d] \le c\,\exp\{- c \,u N^{d-2}\}\,,\label{1.11}
\\[1ex]
&P[\cR_{k_*} \ge u  N^d] \le c\,\exp\{- c \,u N^{d-2}\}\,,\label{1.12}
\end{align}

\medskip\n
for $u > 0$, $N \ge 1$, with $k^* = [c_1 \,u N^{d-2}]$, $k_* = [c_2 \,u N^{d-2}]$, and $c_1 > c_2$. Moreover for $u > 0$, one has:
\begin{equation}\label{1.13}
e^{-cu} \le \liminf\limits_N P[0 \notin X_{[0,uN^d]}] \le \limsup\limits_N P[0 \notin X_{[0,u N^d]}] \le e^{-c^\prime u}\,.
\end{equation}
\end{proposition}

\begin{proof}
We begin with the proof of (\ref{1.11}). As a direct consequence of the invariance principle, we see that for $N \ge 1$, $x \in \wt{B}^c$, $y \in B$,
\begin{equation}\label{1.14}
E_x \Big[\exp\Big\{ - \dis\frac{H_B}{N^2}\Big\}\Big] \le 1 - c, \;E_y \Big[\exp\Big\{ - \dis\frac{T_{\wt{B}}}{N^2}\Big\}\Big] \le 1 - c\,.
\end{equation}

\n
Hence for $k \ge 2$ and $x \in E$, one finds using the strong Markov property and induction:
\begin{equation*}
\begin{split}
E_x \Big[\exp\Big\{ - \dis\frac{\cR_k}{N^2}\Big\}\Big] &  \le (1-c) \,E_x \Big[\exp\Big\{ - \dis\frac{\cD_{k-1}}{N^2}\Big\}\Big] 
\\[1ex]
& \le (1-c)^2 \,E_x \Big[\exp\Big\{ - \dis\frac{\cR_{k-1}}{N^2}\Big\}\Big] \le (1-c)^{2(k-1)}\,.
\end{split}
\end{equation*}

\n
As a result we see that (with the convention below (\ref{1.5}))
\begin{equation}\label{1.15}
E_x \Big[\exp\Big\{ - \dis\frac{\cR_k}{N^2}\Big\}\Big] \le c\,\exp\{- c \,k\}, \;\;\mbox{for $k \ge 0$}\,,
\end{equation}

\n
and therefore for $u > 0$, $N \ge 1$, $k \ge 0$, we find
\begin{equation}\label{1.16}
P[\cR_k \le u N^d] \le c\,\exp\{u N^{d-2} - c\,k\}\,,
\end{equation}
from which (\ref{1.11}) readily follows.

\medskip
We now turn to the proof of (\ref{1.12}). With similar arguments as in the proof of Lemma 1.3 of \cite{DembSzni06}, we see that:
\begin{equation}\label{1.17}
E_x\Big[\exp\Big\{\dis\frac{c}{N^2} \;H_B\Big\}\Big] \le 2, \;E_x \Big[\exp\Big\{\dis\frac{c}{N^2} \;T_{\wt{B}}\Big\}\Big] \le 2,\;\;\mbox{for $N \ge 1$, and $x \in E$}\,.
\end{equation}

\n
Therefore with the strong Markov property and induction we find that for $k \ge 1$:
\begin{equation}\label{1.18}
E\Big[\exp\Big\{\dis\frac{c}{N^2} \;\cR_k\Big\}\Big] \le 2 E  \Big[\exp\Big\{\dis\frac{c}{N} \;\cD_{k-1}\Big\}\Big] \le 4 E  \Big[\exp\Big\{\dis\frac{c}{N^2} \;\cR_{(k-1)}\Big\}\Big] \le 4^k\,.
\end{equation}
It now follows that for $k \ge 0$:
\begin{equation}\label{1.19}
P[\cR_k \ge u N^d] \le  \exp\{- c \,u N^{d-2} + 2(\log 2)k\}\,,
\end{equation}
from which (\ref{1.12}) easily follows.

\medskip
We then prove (1.13). Using a comparison between the Green function of the random walk killed when exiting $\wt{B}$ and of simple random walk in  $\IZ^d$, see for instance (\ref{1.10}), (\ref{1.11}) in Lemma 1.2 of \cite{DembSzni06}, we have
\begin{equation}\label{1.20}
c^\prime(|x- y| + 1)^{-(d-2)} \le P_x [H_y < T_{\wt{B}}] \le c( |x - y| + 1)^{-(d-2)}, \;\;\mbox{for $N \ge 1$, $x,y \in B$}\,.
\end{equation}

\n
Note that for $k \ge 1$, one has:
\begin{equation}\label{1.21}
\{H_0 > \cD_k\} = \{H_0 > \cR_1\} \cap \theta_{\cR_1}^{-1} \{H_0 > T_{\wt{B}}\} \cap \dots \cap \theta^{-1}_{\cR_k} \{H_0 > T_{\wt{B}}\}\,.
\end{equation}

\n
Hence with the strong Markov property and the left hand inequality of (\ref{1.20}) we see that
\begin{equation}\label{1.22}
P[H_0 > \cD_k] \le ( 1- cN^{-(d-2)})^{k}, \;\;\mbox{for $k \ge 1$}\,.
\end{equation}

\n
Similarly we see that for $k \ge 1$, $0 < \e < \frac{1}{8}$,
\begin{equation}\label{1.23}
\begin{split}
P[H_0 > \cD_k] & \ge P\big[X_0 \notin B(0,\e N), \;H_0 > \cD_k\big] 
\\[1ex]
&  \ge \Big( 1 - \dis\frac{c}{N^{d-2}}\Big)^{k-1} P[X_0 \notin B(0, \e N), H_0 > T_{\wt{B}}]
\\
& \ge \big( 1- c N^{-(d-2)}\big)^{k-1} \cdot \big(1 - c(\e N)^{-(d-2)}\big)_+ \Big(1 - \dis\frac{|B(0,\e N)|}{N^d}\Big)\,.
\end{split}
\end{equation}
We can now write for large $N$:
\begin{equation}\label{1.24}
\begin{split}
P[H_0 > u N^d] \le P[\cR_{k_*} \ge u N^d] + P[H_0 > \cD_{k_* - 1}]
\\
\stackrel{(\ref{1.12}), (\ref{1.22})}{\le} c \exp \{ - c \,uN^{d-2}\} + (1 - c N^{-(d-2)})^{k_* -1}\,,
\end{split}
\end{equation}
as well as:
\begin{equation}\label{1.25}
\begin{array}{l}
P[H_0 > u N^d]  \ge P[H_0 > \cD_{k^*}, \cR_{k^*} \ge u N^d] \stackrel{(\ref{1.11}), (\ref{1.23})}{\ge} - c \exp\{-c u\,N^{d-2}\} \;+
\\
\big(1 - c N^{-(d-2)}\big)^{k^*-1} 
\big(1 - c(\e N)^{-(d-2)}\big)  \big(1 - N^{-d} \,|B(0,\e N)|\big) \,.
\end{array}
\end{equation}

\medskip\n
Inserting the value of $k_*$ and $k^*$, see below (\ref{1.12}), we can let $N$ tend to infinity in (\ref{1.24}), (\ref{1.25}) and then $\e$ to $0$ in (\ref{1.25}), and find (\ref{1.13}).
\end{proof}

We now come to the main result of this section that shows the ubiquity of segments of logarithmic size in the vacant set left by the walk at times which are small multiples of $N^d$. As explained in the Introduction, the heuristics underlying this result stems from the coupon collector problem, cf.~below (\ref{0.8}). We refer to (\ref{0.3}) for the definition of the event $\cV_{K,\beta,t}$, with $K > 0$, $0 < \beta < 1$, $t \ge 0$. Our main result is:

\begin{theorem}\label{theo1.2} $(d \ge 4)$

\medskip
For any $K > 0$, $0 < \beta < 1$,
\begin{equation}\label{1.26}
\lim\limits_N \;P[\cV_{K,\beta,u N^d}] = 1, \;\;\mbox{for small $u > 0$}\,.
\end{equation}
\end{theorem}

\begin{proof}
We pick $\beta_1,\beta_2$ such that
\begin{equation}\label{1.27}
0 < \beta_2 < \beta_1 < \beta < 1, \;\mbox{and $2 \beta_1 - \beta_2 < \beta$}\,.
\end{equation}

\n
Using translation invariance and isotropy the claim (\ref{1.26}) follows once we show that
\begin{equation}\label{1.28}
\lim\limits_N N^d \,P\Big[\bigcap\limits_{0 \le m < N^\beta} \{H_{(m + [0,L]) e_1} \le uN^d\}\Big] = 0, \;\mbox{for small  $u > 0$,}
\end{equation}

\n
with the notation $L = [K \log N]$. We now prove (\ref{1.28}), and for this purpose consider the segments $S_i$ in $E$ defined by
\begin{equation}\label{1.29}
S_i = \pi_E \big((2 i [N^{\beta_1 - \beta_2}] + [0,L]) e_1\big), \;\;1 \le i \le \ell, \;\mbox{with} \;\ell \stackrel{\rm def}{=} [N^{\beta_1}] \,,
\end{equation}
and write
\begin{equation}\label{1.30}
S = \bigcup\limits_{1 \le i \le \ell} S_i \,.
\end{equation}

\medskip\n
We want to show that when $u > 0$ is chosen small, with overwhelming probability as $N$ tends to infinity, some of the segments $S_i$, $1 \le i \le \ell$, remain vacant up to time $u N^d$. With the help of (\ref{1.27}), (\ref{1.28}) will then follow. We then introduce, cf.~(\ref{1.10}) for the notation,
\begin{equation}\label{1.31}
\begin{split}
\cS_0  =& \; S, \;\tau_1 = \inf\{k \ge 1; H_S \circ \theta_{\cR_k} < T_{\wt{B}} \circ \theta_{\cR_k}\}\,,
\\[1ex]
\wt{R}_1  = &\;H_S \circ \theta_{\cR_{\tau_1}} + \cR_{\tau_1}, \; \mbox{$j_1 =$ the unique $j \in \{1,\dots,\ell\}$ such that $X_{\wt{\cR}_1} \in S_j$}\,,
\\[1ex]
\cS_1  =& \;S \backslash S_{j_1}, \;\tau_2 = \inf\{k > \tau_1, \;H_{\cS_1} \circ \theta_{\cR_k} < T_{\wt{B}} \circ \theta_{\cR_k}\}\,,
\\[1ex]
\wt{R}_2  = &\;H_{\cS_1} \circ \theta_{\cR_{\tau_2}} + \cR_{\tau_2}, \;\mbox{$j_2 =$ the unique $j \in \{1, \dots,\ell\} \backslash \{j_1\}$ such that} 
\\
&\hspace{-5ex} \mbox{$X_{\wt{R}_2} \in S_j$, and so on until $\cS_{\ell - 1}  = S \;\backslash \bigcup\limits_{j \in \{j_1,\dots,j_{\ell - 1}\}} S_j, \;\tau_\ell ,\;\wt{R}_\ell$ and}
\\[-1ex]
&\hspace{-5ex}\mbox{$j_\ell$ with $\{1,\dots,\ell\} = \{j_1,\dots,j_\ell\}$}\,.
\end{split}
\end{equation}

\medskip\n
In this fashion we label the successive excursions in $B$ and out of $\wt{B}$ giving rise to hits of new segments, and for the time being disregard the fact that possibly more than one segment may be hit during one such excursion. As a straightforward consequence of the above definition, one has
\begin{equation}\label{1.32}
\begin{array}{l}
\wt{R}_i, \; 1 \le i \le \ell, \;\mbox{are $(\cF_n)$-stopping times}\,,
\\[1ex]
\mbox{$j_m, \;1 \le m \le i$, are $\cF_{\wt{R}_i}$-measurable}\,,
\\[1ex]
\cD_{\tau_i} = T_{\wt{B}} \circ \theta_{\wt{R}_i} + \wt{R}_i, \; 1 \le i \le \ell, \;\mbox{are $(\cF_n)$-stopping times as well}\,.
\end{array}
\end{equation}

\medskip\n
With (\ref{1.20}), we see that for large $N$, when $U \subseteq S$
\begin{equation}\label{1.33}
P_x[H_U < T_{\wt{B}}] \le \dis\frac{c_3}{2} \;\dis\frac{|U|}{N^{d-2}}, \;\mbox{for $x \in B \cap \partial (B^c)$}\,,
\end{equation}

\n
and using the reversibility of the walk on $E$,
\begin{equation}\label{1.34}
\begin{split}
P[H_U < T_{\wt{B}}] & \le N^{-d} \dsl_{x \in E, y \in U, k \ge 0} \;P_x [X_k = y, \,T_{\wt{B}} > k]
\\[1ex]
& = N^{-d} \dsl_{x \in E, y \in U, k \ge 0} \;P_y [X_k = x, T_{\wt{B}} > k] = N^{-d} \dsl_{y \in U} \,E_y[T_{\wt{B}}] 
\\[1ex]
& \le \dis\frac{c_3}{2} \;\dis\frac{|U|}{N^{d-2}} \;,
\end{split}
\end{equation}

\medskip\n
using the fact that $\sup_{z \in E} \;E_z [T_{\wt{B}}] \le c N^2$, cf.~the second inequality of (\ref{1.17}), and increasing if necessary the value of $c_3$ in (\ref{1.33}). We now introduce for $U \subseteq S$
\begin{equation}\label{1.35}
\eta(U) = \inf\{k \ge 1; \;H_U \circ \theta_{\cR_k} < T_{\wt{B}} \circ \theta_{\cR_k}\}\,,
\end{equation}

\medskip\n
and note that for $\lambda > 0$, and $2 \le i \le \ell$, as a result of the strong Markov property applied at time $\cD_{\tau_{i-1}}$, cf.~(\ref{1.32}),
\begin{equation}\label{1.36}
E[\exp\{- \lambda(\tau_i - \tau_{i-1})\}\,|\,\cF_{\cD_{\tau_{i-1}}}] = \dsl_{k \ge 1} \,e^{-\lambda k} \;P_{X_{\cD_{\tau_{i-1}}}} [ \eta(\cS_{i-1}) = k]\,,
\end{equation}

\n
where in the last expression $\cS_{i-1}$ is a frozen variable  $(\cF_{\cD_{\tau_{i-1}}}$-measurable).

\medskip
Note that for $z \in \wt{B}^c$, $U \subseteq S$,
\begin{equation}\label{1.37}
\dsl_{k \ge 1} \,e^{-\lambda k} P_z [\eta (U) = k] = (1 - e^{-\lambda}) \,\dsl_{k \ge 1} \,e^{-\lambda k} \,P_z [\eta(U) \le k]\,.
\end{equation}

\n
Moreover for large $N$, with $z$ and $U$ as above
\begin{equation}\label{1.38}
\begin{split}
P_z [\eta(U) \le k]  & = 1 - P_z [\eta(U) > k]  
\\[1ex]
& = 1 - P_z [H_U \circ \theta_{\cR_m} > T_{\wt{B}} \circ \theta_{\cR_m}, \;\mbox{for $1 \le m \le k$}] 
\\
&  \le 1  - \Big(1 - \dis\frac{c_3}{2} \;\dis\frac{|U|}{N^{d-2}}\Big)^k, \;\mbox{for  $k \ge 0$}\,,
\end{split}
\end{equation}

\medskip\n
with the help of (\ref{1.33}) and strong Markov property applied at times $\cR_m$, $1 \le m \le k$. We thus see that under $P_z$, $\eta(U)$ stochastically dominates a geometric variable with success probability
\begin{equation}\label{1.39}
p(U) = c_3 \;\dis\frac{|U|}{N^{d-2}}\;,
\end{equation}

\n
(the factor 2 multiplying the expression subtracted from 1 under the parenthesis in (\ref{1.38}) is there in order to obtain (\ref{1.41}) below). Thus coming back to (\ref{1.37}), we see that for large $N$, $U \subseteq S$, $z \in \wt{B}^c$,
\begin{equation}\label{1.40}
E_z [\exp \{ - \lambda \eta(U)\}] \le \dsl_{k \ge 1} e^{-\lambda k} p(U) \big(1 - p(U)\big)^{k-1} = \dis\frac{e^{-\lambda} p(U)}{1 - (1 - p(U)) e^{-\lambda}}, \;\mbox{for $\lambda \ge 0$}\,.
\end{equation}

\n
In the same fashion, using (\ref{1.33}), (\ref{1.34}), (note that $P[H_S \circ \theta_{\cR_1} < T_{\wt{B}} \circ \theta_{\cR_1}] \le c_3 \,\frac{|S|}{N^{d-2}}$, when bounding the term $m=1$ in the expression corresponding to the second line of (\ref{1.38})),  we find that for large $N$,
\begin{equation}\label{1.41}
E[\exp\{-\lambda \eta  (S)\}] \le \dis\frac{e^{-\lambda} p(S)}{1 - (1 - p(S)) e^{-\lambda}} , \;\;\mbox{for $\lambda \ge 0$}\,.
\end{equation}

\medskip\n
Hence with (\ref{1.36}), (\ref{1.40}), (\ref{1.41}), and the fact that $P$-a.s., for $1 \le i \le \ell$, 
\begin{equation}\label{1.42}
p_i \stackrel{\rm def}{=} p(\cS_{i-1}) \underset{(\ref{1.31})}{\overset{(\ref{1.39})}{=}} \;c_3 (L+1) (\ell - i + 1) \,N^{-(d-2)}\,,
\end{equation}

\n
we find that for large $N$ and $\ell^\prime \le \ell$, $\lambda \ge 0$,
\begin{equation}\label{1.43}
E[\exp\{ - \lambda \tau_{\ell^\prime}\}] \le \prod\limits_{1 \le i \le \ell^\prime} \;\dis\frac{e^{-\lambda} p_i}{1- e^{-\lambda} (1 - p_i)} \;,
\end{equation}

\n
and hence with the notation below (\ref{1.11}),
\begin{equation}\label{1.44}
P[\tau_{\ell^\prime} \le k^*] \le \exp\Big\{\lambda k^* - \lambda \ell^\prime + \dsl_{1 \le i \le \ell^\prime} \;\log \Big(\dis\frac{p_i}{p_i + (1 - e^{-\lambda})(1 - p_i)}\Big)\Big\}\;.
\end{equation}

\medskip\n
We now specify $\ell^\prime$ and $\lambda$ and set, cf.~(\ref{1.27}):
\begin{equation}\label{1.45}
\ell^\prime = \ell - [N^{\beta_2}], \;\;\rho \stackrel{\rm def}{=} d-2 - \dis\frac{(\beta_1 + \beta_2)}{2} \;( > 0) , \;\;\lambda = N^{-\rho}\,.
\end{equation}

\n
As a result for large $N$, with (\ref{1.42}), (\ref{1.27}), we have
\begin{equation}\label{1.46}
\begin{array}{ll}
{\rm i)} & 10^{-3} p_1 \ge N^{-\rho} \ge \dis\frac{(1 - e^{-\lambda})}{2} \ge \dis\frac{N^{-\rho}}{4} \ge 10^3 \,c_3 \,\dis\frac{(L+1)}{N^{d-2}}\;,
\\[2ex]
{\rm ii)} & \fr \ge p_1 \ge p_i, \;\mbox{for $1 \le i \le \ell$}\,,
\\[2ex]
{\rm iii)} & (\ell - \ell^\prime) \;c_3 \;\dis\frac{(L+1)}{N^{d-2}} < \dis\frac{N^{-\rho}}{10} \;.
\end{array}
\end{equation}

\n
As a result we see that for large $N$,
\begin{equation}\label{1.47}
\begin{array}{l}
\dsl_{1 \le i \le \ell^\prime} \;\log \;\Big(\dis\frac{p_i + (1 - e^{-\lambda}) ( 1- p_i)}{p_i}\Big) \ge \dsl_{1 \le i \le \ell^\prime}\; \dis\int_{p_i}^{p_i  + \frac{N^{-\rho}}{4}} \;\dis\frac{dt}{t} \stackrel{(\ref{1.42})}{=}
\\
\\
\dis\int^\infty_0 \;\dsl_{\ell - \ell^\prime < j \le \ell} \;1 \Big\{c_3 \;\dis\frac{(L+1)}{N^{d-2}} \;j < t < c_3 \;\dis\frac{(L+1)}{N^{d-2}} \;j + \dis\frac{N^{-\rho}}{4} \Big\} \;\dis\frac{dt}{t} \;.
\end{array}
\end{equation}

\n
From (\ref{1.46}) it now follows that for large $N$, when $t \in (\frac{N^{-\rho}}{5}$, $c_3 \;\frac{(L+1)\ell}{N^{d-2}})$, the sum inside the above integral is bigger than $c\,N^{-\rho} \frac{N^{d-2}}{L}$. As a result we see that for large $N$, using (\ref{1.27}), the definition of $L$  below (\ref{1.28}), and (\ref{1.45})
\begin{equation}\label{1.48}
\begin{array}{l}
\dsl_{1 \le i \le \ell^\prime} \;\log \;\Big(\dis\frac{p_i + (1 - e^{-\lambda}) (1- p_i)}{p_i}\Big) \ge c \,N^{-\rho} \;\dis\frac{N^{d-2}}{L} \;\log (5c_3 (L+1) \,\ell N^{\rho - (d-2)}) 
\\[2ex]
\ge \dis\frac{c}{K} \;N^{d-2-\rho} \,\big(\beta_1 + \rho - (d-2)\big) = \frac{c}{K} \;(\beta_1 - \beta_2) \;N^{\frac{\beta_1 + \beta_2}{2}}\;.
\end{array}
\end{equation}

\n
Inserting this bound in (\ref{1.44}), with (\ref{1.45}) we see that for large $N$:
\begin{equation}\label{1.49}
P[\tau_{\ell^\prime} \le k^*] \le \exp\Big\{ c_1 \,u N^{\frac{\beta_1 + \beta_2}{2}} - \dis\frac{c}{K} \;(\beta_1 - \beta_2) \;N^{\frac{\beta_1 + \beta_2}{2}}\Big\}\;.
\end{equation}

\n
We will now take care of the issue (cf.~below (\ref{1.31})), of additional hits during the time intervals $[\wt{R}_i,\cD_{\tau_i}]$, $1 \le i \le \ell^\prime$, of non-previously hit segments. With this objective in mind we thus define for $1 \le i \le \ell^\prime$,
\begin{equation}\label{1.50}
N_i = \dsl_{j \notin \{j_1,\dots,j_i\}} 1\{H_{S_j} \circ \theta_{\wt{R}_i} < T_{\wt{B}} \circ \theta_{\wt{R}_i}\} \,,
\end{equation}

\n
so that $N_i$ is $\cF_{\cD_{\tau_i}}$-measurable, cf.~(\ref{1.32}). Note that $P$-a.s, $X_{\wt{R}_i} \in S_{j_i}$, for $1 \le i \le \ell$, cf.~(\ref{1.31}), and with the strong Markov property at time $\wt{R}_i$, one finds
\begin{equation}\label{1.51}
\begin{split}
P[N_i \ge m \,|\,\cF_{\wt{R}_i}] & \le P_{X_{\wt{R}_i}} \Big[\dsl_{1 \le j \le \ell} \,1\{H_{S_j} < T_{\wt{B}}\} \ge 1 + m\Big]
\\[1ex]
& \le P_{X_{\wt{R}_i}} [V_m < T_{\wt{B}}] \,,
\end{split}
\end{equation}

\n
where $V_m, m \ge 1$, denote the successive times of visit of the walk to distinct segments $S_j$, $1 \le j \le \ell$.

\medskip
It now follows from (\ref{1.20}), (\ref{1.29}) that for large  $N$, when $z \in S$.
\begin{equation}\label{1.52}
P_z [V_1 < T_{\wt{B}}] \le c(L+1) \;\dsl_{k \ge 1} \;(k\,N^{\beta_1 - \beta_2})^{-(d-2)} \le c_4 \,L N^{-(d-2)(\beta_1 - \beta_2)} \stackrel{\rm def}{=} \;p\,.
\end{equation}

\n
Coming back to (\ref{1.51}) we thus see with the strong Markov property that when $N$ is large,
\begin{equation}\label{1.53}
P[N_i \ge m \,| \,\cF_{\wt{R}_i}] \le p^m,\;\mbox{for $m \ge 0$} \,,
\end{equation}

\medskip\n
i.e. conditionally on $\cF_{\wt{R}_i}$, $N_i$ is stochastically dominated by a modified geometric distribution with success parameter $p$. Therefore when $\lambda^\prime$ is such that, cf.~(\ref{1.52}), $e^{\lambda^\prime} \,p < 1$, we find that for large $N$, $\lambda^\prime$ as above and  $1 \le i \le \ell^\prime$: 
\begin{equation}\label{1.54}
E\big[\exp\{\lambda^\prime N_i\} \,| \,\cF_{\wt{R}_i}\big] \le \dsl_{m \ge 0} \;(1-p) \,p^m \,e^{\lambda^\prime m} = \dis\frac{1- p}{1 - e^{\lambda^\prime} p} \;,
\end{equation}

\n
so that using induction and $\cF_{\cD_{\tau_{i-1}}} \subseteq \cF_{\wt{R}_i}$ for $2 \le i \le \ell^\prime$, we see in view of the measurability of $N_i$ asserted below (\ref{1.50}) that for large $N$:
\begin{equation}\label{1.55}
E\Big[ \exp\Big\{ \dsl_{1 \le i \le \ell^\prime} N_i\Big\}\Big] \le \Big(1 + \dis\frac{(e - 1)\,p}{1 -  ep}\Big)^{\ell^\prime}\stackrel{(\ref{1.52})}{\le}  \exp \{ c\,\ell^\prime p\}\,,
\end{equation}

\medskip\n
and hence with (\ref{1.45}), the value of $p$ in (\ref{1.52}), and the fact that $d \ge 4$:
\begin{equation}\label{1.56}
\begin{split}
P\Big[\dsl_{1 \le i \le \ell^\prime} \,N_i \ge \mbox{\f $\dis\frac{\ell - \ell^\prime}{2}$}\Big] & \le \exp\Big\{ - \mbox{\f $\dis\frac{N^{\beta_2}}{4}$} + c \,\ell^\prime p\Big\}
\\
& \le \exp\Big\{ - \mbox{\f $\dis\frac{N^{\beta_2}}{4}$} + c \,N^{\beta_2} L\,N^{-(d-3)(\beta_1 - \beta_2)}\Big\}   \le \exp \Big\{ - \mbox{\f $\dis\frac{N}{8}^{\beta_2}$}\Big\} \,.
\end{split}
\end{equation}

\n
To conclude the proof of (\ref{1.28}), we observe that with (\ref{1.31}) for large $N$, on the event $\{\tau_{i-1} < k < \tau_i\}$, where $1 \le i \le \ell^\prime$ and $\tau_0 = 0$ by convention, $X_n \notin \cS_{i-1}$ for $\cR_k \le n < \cR_{k+1}$. As a result on the event $\{\cR_{k^*} > u N^d\} \cap \{\tau_{\ell^\prime} > k^*\} \cap \{\sum_{1 \le i \le \ell^\prime} N_i < \frac{\ell - \ell^\prime}{2}\}$ at least $[\frac{\ell - \ell^\prime}{2}]$ segments $S_i$, $1 \le i \le \ell$, have not been visited by the walk up to time $u N^d$, so that when $N$ is large the above event lies in the complement of the event that appears in (\ref{1.28}). Collecting the bounds (\ref{1.11}), (\ref{1.49}), (\ref{1.56}), we obtain (\ref{1.28}). As already explained, this yields our claim (\ref{1.26}), so that Theorem \ref{theo1.2} is now proved.
\end{proof}

\begin{remark}\label{rem1.3} \rm Concerning the large $u$ regime, let us point out that when $d \ge 4$, given $K > 0$, the vacant set left by the walk at time $uN^d$, typically for large $N$ does not contain any segment of length $[K \log N]$, if $u$ is chosen large enough. Indeed when $8 L \le N$, and $U \subseteq B$ is the segment $U = [0,L]\,e$, with $|e| = 1$,
\begin{equation}\label{1.57}
\begin{split}
P_x[H_U < T_{\wt{B}}] & \ge E_x \Big[\dsl_{n=0}^{T_{\wt{B}} - 1} 1 \{X_n \in U\}\Big] \Big/
\sup\limits_{y \in U} \,E_y \Big[\dsl_{n=0}^{T_{\wt{B}} - 1}1 \{X_n \in U\}\Big]
\\
& \ge c\;\dis\frac{L}{N^{d-2}} , \;\mbox{for any $x \in B$}\,,
\end{split}
\end{equation}

\medskip\n
using the strong Markov property at time $H_U$ for the first inequality, and for the second inequality bounds in the Green function of the walk killed when exiting $\wt{B}$, see (\ref{1.20}) and also (\ref{1.11}) of \cite{DembSzni06}. With a straightforward modification of (\ref{1.24}), we thus find that when $8 L \le N$, and $u > 0$,
\begin{equation}\label{1.58}
\begin{split}
P[H_U > uN^d] & \le c \exp\{- c u N^{d-2}\} + \Big(1 - c\; \dis\frac{L}{N^{d-2}}\Big)_+^{k_* - 1}
\\
& \le c \exp\{- c u L\}\,,
\end{split}
\end{equation}

\n
using the value of $k_*$ below (\ref{1.12}). Hence choosing $L = [ c_*\, \frac{\log N}{u}]$, with $c_*$ a large enough constant, we see that for $u > 0$,
\begin{equation}\label{1.59}
\mbox{$\lim\limits_N P\Big[X^c_{[0,uN^d]}$ contains some segment of length $\Big[c_* \,\dis\frac{\log N}{u}\Big]\Big] = 0$}\,.
\end{equation}

\medskip\n
So for large $N$ the vacant set typically does not contain segments of length $K \log N$, if $u$ is chosen large enough. \hfill $\square$
\end{remark}

The theorem we have just proved will enter as a step when showing in the next section that the giant component we define,  with overwhelming probability occurs in the regime of parameters we consider.

 \section{Exponential bound and giant component}
 \setcounter{equation}{0}
 
We derive in this section an exponential bound on the probability that the walk covers certain subsets of $E$ by times  that are small multiple of $N^d$, cf.~Theorem \ref{theo2.1}. This bound plays an important role in the construction of the giant component typically present in the vacant set left by the walk at such times. We also refer to Remark \ref{rem2.4} where it is explained why some restrictions are needed on the class of sets to which the exponential bound applies.

\medskip
We recall (\ref{1.2}) for the definition of $\cL_m$, $1 \le m \le d$, and define for $1 \le m \le d$,
\begin{equation}\label{2.1}
\begin{split}
\cA_m = &\; \mbox{the collection of non-empty subsets $A$ of $E$ such that $A \subseteq F$,}\\
&\; \mbox{for some $F \in \cL_m$}\,.
\end{split}
\end{equation}

\medskip\n
Clearly $\cA_m$ increases with $m$, and $\cA_d$ is the collection of non-empty subsets of $E$. We will especially be interested in $\cA_2$. We also recall the notation $q(\nu)$ in (\ref{0.10}) and below (\ref{1.8}). The next theorem contains the key exponential estimate.

\begin{theorem}\label{theo2.1} $(d \ge 4$, $1 \le m \le d-3)$

\medskip
When $\lambda > 0$ is such that
\begin{align}
&\chi \stackrel{\rm def}{=} \;e^{2 \lambda} \Big(\mbox{\f $\dis\frac{m}{d}$} + \Big(1 - \mbox{\f $\dis\frac{m}{d}$}\Big) \,q(d-m) \Big) < 1, \;\mbox{then for $u > 0$,} \label{2.2}
\\[1ex]
&\limsup\limits_N \;\sup\limits_{A \in \cA_m} |A|^{-1} \log \Big( E \Big[ \exp\Big\{\lambda \,\dsl_{x \in A} \;1\{H_x \le u N^d\big\} \Big\} \Big]\Big) \le c_5 \,u \;\mbox{\f $\dis\frac{e^{2\lambda}-1}{1-\chi}$}\;, \label{2.3}
\end{align}
and there exist $N_1(d,m,\lambda) \ge 1$, $u_1(d,m,\lambda) > 0$, such that for $N \ge N_1$:
\begin{equation}\label{2.4}
P[X_{[0,u_1 N^d]} \supseteq A] \le \exp\{ - \lambda \,|A|\}, \;\mbox{for all $A \in \cA_m$}\;.
\end{equation} 
\end{theorem}

We refer to Remark \ref{rem2.4} below for an explanation on why some restriction on the class of subsets $A$ that appear in (\ref{2.4}) is needed.

\begin{proof}
We begin with the proof of (\ref{2.3}). We consider $N \ge 1$, $u > 0$, $A \in \cA_m$, $1 \le m \le d-3$. Roughly speaking we chop the time interval $[0,[u N^d]]$ into successive intervals of length $N^2$, except maybe for the last one, and write for $\lambda > 0$:
\begin{equation}\label{2.5}
\begin{array}{l}
E\Big[\exp\Big\{\lambda \dsl_{x \in A} 1\{H_x  \le uN^d\}\Big\}\Big] \le 
\\[1ex]
E \Big[\exp\Big\{\lambda \dsl_{k N^2 \le u N^d}\;\; \dsl_{x \in A} \,1\{H_x < N^2\} \circ \theta_{kN^2}\Big\}\Big] \le \sqrt{a_1} \;\sqrt{a_2}, \;\mbox{where}
 \end{array}
\end{equation}
\begin{equation}\label{2.6}
\begin{array}{l}
a_1 \stackrel{\rm def}{=} E\Big[\exp \Big\{ 2 \lambda \dsl_{k \;{\rm even}, \,kN^2 \le uN^d} \;\dsl_{x \in A} \,1\{H_x < N^2\} \circ \theta_{k N^2}\Big\}\Big]\,, 
\\[1ex]
a_2\stackrel{\rm def}{=} E\Big[\exp \Big\{ 2 \lambda \dsl_{k \;{\rm odd}, \,kN^2 \le uN^d} \;\dsl_{x \in A} \,1\{H_x < N^2\} \circ \theta_{k N^2}\Big\}\Big]\,.
\end{array}
\end{equation}

\n
We first bound $a_1$. To this end we define
\begin{align}
& k_0 = \max\{k \ge 0, \;2 k \,N^2 \le u N^d\}, \;\mbox{and} \label{2.7}
\\[1ex]
& \phi(z) = E_z \Big[\exp\Big\{ 2 \lambda  \,\dsl_{x \in A} \,1\{H_x < N^2\}\Big\}\Big] \;(\ge 1), \;\mbox{for $z \in E$}\,. \label{2.8}
\end{align}

\n
Applying the strong Markov property at time $H_A$, we find:
\begin{equation}\label{2.9}
\phi(z) \le P_z [H_A \ge N^2] + E_z [H_A < N^2, \;\phi(X_{H_A})], \;\mbox{for $z \in E$}\,.
\end{equation}

\n
With the simple Markov property applied at time $2k_0  \,N^2$ and then at time $(2k_0 - 1) \,N^2$, we see that when $k_0 \ge 1$:
\begin{equation}\label{2.10}
a_1 = E\Big[ \exp\Big\{ 2\lambda \,\dsl_{0 \le k < k_0} \;\dsl_{x \in A} \,1\{H_x < N^2\} \circ \theta_{2k N^2}\Big\} \,E_{X_{(2 k_0 -1) N^2}} [\phi(X_{N^2})]\Big]\,.
\end{equation}

\n
Note that for $z \in E$, one has:
\begin{equation}\label{2.11}
\begin{array}{lcl}
E_z [\phi(X_{N^2})] &\hs \stackrel{(\ref{1.6})}{=}&\hs \dsl_{y \in E} \;p_{N^2} (z,y) \,\phi(y) 
\\[1ex]
& \hs \stackrel{(\ref{2.9})}{\le} &\hs 1 + \dsl_{y \in E} \,p_{N^2}(z,y) \,E_y[H_A < N^2; \phi(X_{H_A}) - 1]
\\[3ex]
&\hs \le & \hs 1 + \dis\frac{c}{N^d} \;\dsl_{y \in E} \;\dsl_{0 \le k < N^2} \;P_y [X_k \in A] \;(\|\phi\|_\infty - 1)
\\[1ex]
&\hs \le &\hs 1 + \dis\frac{c}{N^{d-2}} \;|A| \;(\|\phi\|_\infty - 1) \le \exp\Big\{c \;\dis\frac{|A|}{N^{d-2}} \;(\|\phi\|_\infty - 1)\Big\}\,,
\end{array}
\end{equation}

\medskip\n
where in the third line we used that
\begin{equation}\label{2.12}
\sup\limits_{x,y \in E} \,N^d \,p_{N^2} (x,y) \le c\,,
\end{equation}

\n
as follows from standard upper bounds on the transition density of simple random walk on $\Z^d$, cf.~(2.4) of \cite{GrigTelc01}. With an even simpler (and similar) argument we also have
\begin{equation}\label{2.12a}
E[\phi(X_0)] \le \exp\Big\{c\;\dis\frac{|A|}{N^{d-2}} \;(\|\phi\|_\infty - 1)\Big\}\,.
\end{equation}

\n
Therefore using induction together with (\ref{2.11}), and (\ref{2.12a}) to handle the term corresponding to $k = 0$ in (\ref{2.10}), we see that
\begin{equation}\label{2.13}
a_1 \le \exp\Big\{ (k_0 + 1) \;c\,\dis\frac{|A|}{N^{d-2}} \;(\|\phi\|_\infty - 1)\Big\}\,.
\end{equation}

\medskip\n
A similar bound holds for $a_2$, and with (\ref{2.5}) we thus find
\begin{equation}\label{2.14}
E\Big[\exp\Big\{ \lambda \,\dsl_{x \in A} \, 1\{H_x \le u N^d\}\Big\}\Big] \le \exp\Big\{c (u N^{d-2} + 1) \;\dis\frac{|A|}{N^{d-2}} \;(\|\phi\|_\infty - 1)\Big\}\,.
\end{equation}

\medskip\n
We will now seek an upper bound on $\|\phi\|_\infty$.

\begin{lemma}\label{lem2.2} $\big(d \ge 4, 1 \le m \le d-3, e^{2\lambda} \;\mbox{\f $\dis\frac{m}{d}$} < 1, N \ge  2\big)$
\begin{equation}\label{2.15}
\|\phi\|_\infty \le \dis\frac{e^{2\lambda}}{1 - e^{2 \lambda} \;\frac{m}{d}} \;\Big(1 - \mbox{\f $\dis\frac{m}{d}$}\Big) \;\big(1 + q_N \,(\|\phi\|_\infty - 1)\big) \,,
\end{equation}

\medskip\n
where (with hopefully obvious notations),
\begin{equation}\label{2.16}
q_N \stackrel{\rm def}{=} P_{e_1}^{(\IZ/N\IZ)^{d-m}} [H_0 < N^2]\,.
\end{equation}
\end{lemma}

\begin{proof}
Consider $F \in \cL_m$ such that $A \subseteq F$, and introduce, cf. (\ref{1.3}),
\begin{equation}\label{2.17}
R_F \stackrel{\rm def}{=} H_F \circ \theta_{T_F} + T_F\,,
\end{equation}

\n
the return time to $F$. Since $A \subseteq F$, for $z \in E$, we find:
\begin{equation}\label{2.18}
\begin{split}
&\phi(z)  \le E_z \Big[\exp\Big\{2 \lambda \Big(T_F + 1\{R_F < N^2\} \Big(\dsl_{x \in A} \,1\{H_x < N^2\} \circ \theta_{R_F}\Big)\Big)\Big\}\Big]
\\[1ex]
&=  E_z \Big[ \exp \{2 \lambda T_F\} \Big(1{\{R_F \ge N^2\}} + 1{\{R_F < N^2\}} \,\exp\Big\{2 \lambda \,\dsl_{x \in A} \,1 \{H_x < N^2\} \circ \theta_{R_F}\Big\}\Big)\Big]
\\[1ex]
&=  E_z \big[ \exp \{2 \lambda T_F\} \big(1 +  1\{R_F <  N^2\} \big(\phi(X_{R_F}) - 1\big)\big)\big]
\\[2ex]
&\le  E_z \big[ \exp \{2 \lambda T_F\}\big] + E_z \big[\exp \{2 \lambda T_F\} \,P_{X_{T_F}} [H_F < N^2]\big] (\|\phi\|_\infty - 1)\,,
\end{split}
\end{equation}

\medskip\n
where we used the strong Markov property at time $R_F$ in the third line. Considering the motion of $X$ in the directions ``transversal to $F$'', we have:
\begin{equation}\label{2.19}
\mbox{for $z \in E$, $P_z$-a.s., $P_{X_{T_F}} [H_F < N^2] \le q_N$}\,.
\end{equation}

\n
When $z \in F$, $T_F$ has geometric distribution with success probability $1 - \frac{m}{d}$, so that for $\lambda$ as indicated above
\begin{equation}\label{2.20}
\begin{array}{l}
E_z [\exp\{2 \lambda T_F\}] = 
\\[1ex]
\dsl_{k \ge 1} \;\Big(1 - \mbox{\f $\dis\frac{m}{d}$}\Big)\Big(\mbox{\f $\dis\frac{m}{d}$}\Big)^{k-1} \,e^{2 \lambda k} = e^{2 \lambda} \Big(1 - \mbox{\f $\dis\frac{m}{d}$}\Big) \Big(1 - e^{2 \lambda} \;\mbox{\f $\dis\frac{m}{d}$}\Big)^{-1}, \;\mbox{for $z \in F$},
\end{array}
\end{equation}

\medskip\n
whereas $T_F = 0$, $P_z$-a.s., when $z \notin F$. Hence coming back to the last line of (\ref{2.18}), we obtain (\ref{2.15}).
\end{proof}

In the next lemma we relate $q_N$ of (\ref{2.16}) to $q(d-m)$, cf.~(\ref{0.10}).
\begin{lemma}\label{lem2.3}
\begin{equation}\label{2.21}
\limsup\limits_N \;q_N \le q(d-m)\,.
\end{equation}
\end{lemma}

\begin{proof}
We denote with $W$ the discrete cube image of $V \stackrel{\rm def}{=} [-\frac{N}{4}, \frac{N}{4}]^{d-m} \cap \IZ^{d-m}$ under the canonical projection onto $(\IZ / N\IZ)^{d-m}$. We have
\begin{equation}\label{2.22}
\begin{split}
q_N &\le  P_{e_1}^{(\IZ/N\IZ)^{d-m}}[H_0 < T_W] + E_{e_1}^{(\IZ/N\IZ)^{d-m}} \big[P^{(\IZ/N\IZ)^{d-m}}_{X_{T_W}} [H_0 < N^2]\big]
\\[1ex]
& \le q(d-m) + \sup\limits_{z \in \partial W} \;P_z^{(\IZ/N\IZ)^{d-m}} [H_0 < N^2]\,.
\end{split}
\end{equation}

\medskip\n
One has the classical upper bound, cf.~for instance (\ref{2.4}) of \cite{GrigTelc01},
\begin{equation}\label{2.23}
P_x^{\IZ^{d-m}}[X_k = y] \le \dis\frac{c(m)}{k^{\frac{d-m}{2}}} \;\exp\Big\{-c(m) \;\dis\frac{|y -x |^2}{k}\Big\}, \;\mbox{for  $k \ge 1$, $x,y \in \IZ^{d-m}$}\,,
\end{equation}

\medskip\n
(using the convention concerning constants stated at the end of the Introduction). Hence for large $N$ we obtain
\begin{equation}\label{2.24}
\begin{array}{l}
\sup\limits_{z \in \partial W} \;P_z^{(\IZ / N\IZ)^{d-m}} [H_0 < N^2] = \sup\limits_{z \in \partial V} \;P_z^{\IZ^{d-m}} [H_{N\IZ^{d-m}} < N^2] \stackrel{(\ref{2.23})}{\le}
\\[2ex]
\sup\limits_{z \in \partial V} \;\dsl_{y \in \IZ^{d-m}} \;\dsl_{2 \le k < N^2} \;\dis\frac{c(m)}{k^{\frac{d-m}{2}}} \;\exp\Big\{ - c(m) \;\dis\frac{|N y-z|^2}{k} \Big\} \le 
\\[3ex]
\sup\limits_{z \in \partial V}  \;\dsl_{y \in \IZ^{d-m}} \;\dis\int_{0}^{N^2} \; \;\dis\frac{c(m)}{s^{\frac{d-m}{2}}} \;\exp\Big\{ - c(m) \;\dis\frac{|N y-z|^2}{s} \Big\} \le
\\[3ex]
\sup\limits_{w \in \frac{\partial V}{N}} \;c(m)\,N^{-(d-m-2)} \,\dis\int^1_0 \;\dsl_{y \in \IZ^{d-m}} \;t^{- (\frac{d-m}{2})} \;\exp\Big\{-c(m) \;\dis\frac{|y - w|^2}{t} \Big\} \;dt\,.
\end{array}
\end{equation}

\medskip\n
We can now split the sum under the integral, keeping on one hand $y \in \IZ^{d-m}$ with $|y| \ge c(m)$, so that
\begin{equation*}
\mbox{$|y - w|^2 \ge c(m) \,|y^\prime|^2$,  for $y^\prime \in y + [0,1]^{d-m}$ and $w \in \dis\frac{\partial V}{N} \;( \subseteq [-1,1]^{d-m})$}\,,
\end{equation*}

\n
and hence for $t \in (0,1]$, $w \in \frac{\partial V}{N}$:
\begin{equation}\label{2.25}
\begin{array}{l}
\dsl_{|y| \ge c(m)} t^{-\frac{(d-m)}{2}} \exp\Big\{- c(m) \;\dis\frac{|y - w|^2}{t}\Big\} \le
\\[1ex]
\dsl_{y \in \IZ^{d-m}} \;\dis\int_{y + [0,1]^{d-m}} \;t^{-\frac{(d-m)}{2}} \;\exp\Big\{ - c(m) \;\dis\frac{|y^\prime|^2}{t}\Big\} \,dy^\prime \le c(m)
\end{array}
\end{equation}

\medskip\n
and on the other hand we consider the finitely many terms corresponding to $|y| < c(m)$. For these terms we also have in view of the definition of $V$:
\begin{equation*}
\inf \Big\{|y - w|^2, \;w \in \mbox{\f $\dis\frac{\partial V}{N}$}, \;y \in \IZ^{d-m}\Big\} \ge c(m) > 0\,,
\end{equation*}
so that for $w \in \frac{\partial V}{N}$,
\begin{equation}\label{2.26}
\dis\int^1_0 \;\dsl_{|y| < c(m)} \;t^{-\frac{(d-m)}{2}} \;\exp\Big\{ - c(m) \;\dis\frac{|y - w|^2}{t} \Big\} \;dt \le c(m) \,.
\end{equation}

\medskip\n
Thus coming back to the last line of (\ref{2.24}) we find for large $N$
\begin{equation}\label{2.27}
\sup\limits_{z \in \partial W} \;P_z^{(\IZ / N \IZ)^{d-m}} [H_0 < N^2] \le c(m) \,N^{-(d-m-2)} \,,
\end{equation}

\n
and since $d-m-2 \ge 1$ by assumption, letting $N$ tend to infinity in (\ref{2.22}) we find (\ref{2.21}).
\end{proof}

With (\ref{2.15}), (\ref{2.21}), it follows with a straightforward computation that when $\lambda > 0$ satisfies (\ref{2.2}):
\begin{equation}\label{2.28}
\limsup\limits_N \;\sup\limits_{A \in \cA_m} \,(\| \phi \|_\infty - 1) \le \dis\frac{e^{2 \lambda} - 1}{1 - \chi} \;.
\end{equation}

\n
Coming back to (\ref{2.14}), taking logarithms and dividing by $|A|$, the claim (\ref{2.3}) readily follows.

\medskip
We now turn to the proof of (\ref{2.4}). We pick $\wt{\lambda}(d,m,\lambda) > \lambda$ and $\wt{q} (d,m,\lambda) > q$ so that
\begin{equation}\label{2.29}
1 - e^{2 \wt{\lambda}} \Big(\dis\frac{m}{d} + \Big( 1 - \dis\frac{m}{d}\Big) \;\wt{q}\Big) = \fr \;(1 - \chi) \,.
\end{equation}

\n 
Applying (\ref{2.3}) with $\wt{\lambda}$ (for which (\ref{2.2}) holds) we see that for $u > 0$, $N \ge N_2 (d,m,\lambda,u)$ and any $A \in \cA_m$
\begin{equation}\label{2.30}
P[X_{[0,uN^d]} \supseteq A] \le \exp\Big\{ - \wt{\lambda} \,|A| + c u \;\dis\frac{e^{2 \wt{\lambda}} -1}{1- e^{2 \wt{\lambda}} (\frac{m}{d} + (1 - \frac{m}{d}) \,\wt{q})}\;|A|\Big\}\,.
\end{equation}

\medskip\n
Choosing $u = u_1(d,m,\lambda)$ small enough, and setting $N_1(d,m,\lambda) = N_2(d,m,\lambda,u_1)$, we obtain (\ref{2.4}).
\end{proof}

\begin{remark}\label{rem2.4}\rm  ~

\medskip\n
1) Let us mention that it is straightforward to argue in Lemma \ref{lem2.3} that $\liminf_N q_N \ge q(d-m)$, so that (\ref{2.22}) can be sharpened into
\begin{equation}\label{2.22'}
\lim\limits_N q_N = q(d-m)\,,
\end{equation}
although we do not use this sharpened limiting result here.

\bigskip\n
2) As we now explain there is no exponential bound of type (\ref{2.4}) valid uniformly for all $A \in \cA_d$ (i.e. all non-empty subsets of $E$), when $N$ is large, no matter how small $\lambda > 0$ is chosen. Indeed when $\rho \in (0,1)$ and $A_L = \pi_E ([-L,L]^d)$, with $L = [N^\rho]$, a qualitatively similar calculation as in Proposition 2.7 in Chapter 3 of \cite{Szni98a}, see in particular p.~114, (the calculation in \cite{Szni98a} is performed in a Brownian motion setting), shows that for large $N$, $T = [c \,L^d \log L]$, and for all $x \in A_L$, 
\begin{equation}\label{2.31}
\begin{array}{l}
P_x [X_{[0,T]} \supseteq A_L] \ge P_x [X_{[0,T]} \supseteq A_L, \;T_{A_{2L}} > T] \ge \fr \;P_x [T_{A_{2L}} > T]
\\[1ex]
\ge c \,\exp \Big\{- \dis\frac{c}{L^2}\;T\Big\} \ge c\,\exp\Big\{- c \;\dis\frac{|A_L|}{L^2} \;\log L\Big\}\,.
\end{array}
\end{equation}

\medskip\n
Moreover with standard transition density estimates, cf.~(\ref{2.4}) of \cite{GrigTelc01}, one has:
\begin{equation*}
\inf\limits_{z \in E} \;P_z [H_{A_L} < N^2] \ge c\,\Big(\dis\frac{L}{N}\Big)^{d-2} \,,
\end{equation*}

\medskip\n
so that using Markov property at times $k N^2$, one finds for large $N$: 
\begin{equation}\label{2.32}
P\Big[H_{A_L} > \dis\frac{u}{2} \;N^d\Big] \le \Big(1 - c \Big(\dis\frac{L}{N}\Big)^{d-2}\Big)^{[\frac{u}{2} \,N^{d-2}]} \le \fr \,.
\end{equation}

\medskip\n
As a result we see that for any $u > 0$ and $0 < \rho < 1$,
\begin{equation}\label{2.33}
\liminf\limits_N \;\big(|A_L|^{\frac{d-2}{d}} \;\log \,|A_L|)^{-1} \,\log P[X_{[0,uN^d]} \supseteq A_L] > - \infty\,,
\end{equation}
and hence
\begin{equation}\label{2.34}
\lim\limits_N  \;\sup\limits_{A \in \cA_d} \;|A|^{-1} \log P[X_{[0,uN^d]} \supseteq A] = 0\,.
\end{equation}

\medskip\n
This explains why some restriction on the class of subsets $A$ entering (\ref{2.4}) is needed. \hfill $\square$
\end{remark}

We now turn to the applications of Theorem \ref{theo2.1} to the construction of the giant component in the vacant set left by the walk at times that are small multiples of $N^d$. We recall that $\star$-nearest neighbor paths have been defined at the beginning of Section 1, and write:
\begin{align}
a(n) = &\; \mbox{the cardinality of the collection of  $\star$-nearest neighbor self-avoiding   } \label{2.35}
\\[-0.5ex]
&\; \mbox{paths on $\IZ^2$, starting at the origin, with $n$ steps}\,.\nonumber
\end{align}

\n
One has the easy upper bound:
\begin{equation}\label{2.36}
a(n) \le 8.7^{n-1}, \;\mbox{for $n \ge 1$}\,.
\end{equation}

\n
We now define for $N \ge 1$, $K > 0$, $t \ge 0$, the event, cf.~(\ref{1.2}) for the notation,
\begin{align}
\cU_{K,t} = \{ &\mbox{for any $F \in \cL_2$, and connected subsets $O_1,O_2$ of $F \backslash X_{[0,t]}$}, \label{2.37}
\\
& \mbox{with $|\cdot|_\infty$-diameter at least $[K \log N]$, $O_1$ and $O_2$ are in} \nonumber 
\\
&\mbox{the same component of $F \backslash X_{[0,t]}\}$}\,.\nonumber
\end{align}

\n
The above event will be useful in singling out the giant component. The next event will be convenient in the derivation of lower bounds on the relative volume of the giant component in Section 4. For $N \ge 1$, $K > 0$, $x \in E$, $t \ge 0$, we define with the notation of the beginning of Section 1:
\begin{align}
\cC_{K,x,t} = \{ &\mbox{for some $F \in \cL_2$, with $x \in F$, there is a nearest} \label{2.38}
\\
& \mbox{neighbor path in $F\backslash X_{[0,t]}$, from $x$ to $S(x,[K \log N])\}$} \,. \nonumber 
\end{align}

\n
We can now state
\begin{corollary}\label{cor2.5}
There is a smallest $d_0 \ge 5$, such that
\begin{equation}\label{2.39}
\mu \stackrel{\rm def}{=} 49 \Big(\mbox{\f $\dis\frac{2}{d}$} + \Big(1 - \mbox{\f $\dis\frac{2}{d}$}\Big) \,q(d-2)\Big) < 1, \;\mbox{for $d \ge d_0$} \,.
\end{equation}
For $d \ge d_0$, there is a constant $c_0 > 0$, cf.~(\ref{2.45}), such that
\begin{align}
&\lim\limits_N \;P[\cU_{c_0,uN^d}] = 1, \;\mbox{for small $u > 0$, and} \label{2.40}
\\[1ex]
&\lim\limits_{u \rightarrow 0} \;\liminf\limits_N \;P[\cC_{c_0,0,uN^d}] = 1\,, \label{2.41}
\end{align}

\n
(and of course $P[\cC_{c_0,x,uN^d}]  = P[\cC_{c_0,0,uN^d}]$ for all $x \in E$).
\end{corollary}

\begin{proof}
One knows, cf.~(5.4) in \cite{Mont56}, that $q(\cdot)$ has the asymptotic behavior:
\begin{equation}\label{2.42}
q(\nu) \sim (2\nu)^{-1}, \;\mbox{as $\nu \rightarrow \infty$}\,,
\end{equation}

\n
so that (\ref{2.39}) straightforwardly follows. Now consider $d \ge d_0$, and choose $\lambda_0(d)$ such that
\begin{equation}\label{2.43}
e^{\lambda_0} \stackrel{\rm def}{=} 7 \mu^{-\frac{1}{4}} ( > 7), \;\mbox{so that $e^{2 \lambda_0} \,\Big(\mbox{\f $\dis\frac{2}{d}$} + \Big(1 - \mbox{\f $\dis\frac{2}{d}$}\Big) \;q(d-2)\Big) < 1$} \,.
\end{equation}

\n
When N is large, on $\cU^c_{K,uN^d}$, one can find $F \in \cL_2$ and $O_1,O_2 \subseteq F \backslash X_{[0,uN^d]}$, distinct connected components of $F \backslash X_{[0,uN^d]}$ with $|\cdot|_\infty$-diameter at least  $[K\log N]$. We can then introduce $\wh{O_i}, i=1,2$, the inverse images of $O_i$ under an ``affine projection'' of $\IZ^2$ onto $F$. Considering separately the case when at least one of the $\wh{O_i}, i=1,2$, has bounded components, (necessarily of $|\cdot|_\infty$-diameter at least  $[K\log N]$), or both of the $\wh{O_i}$ have unbounded components, one can construct a $\star$-nearest neighbor self-avoiding path $\pi$ with $[K\log N]$ steps in $\partial O_1 \cap F$ or $\partial O_2 \cap F$ ( $\subseteq F \cap X_{[0,uN^d]}$), see also Proposition 2.1, p.387, in  \cite{Kest82}. Therefore for $u < u_0 = u_1 (d,m=2,\lambda = \lambda_0)$, cf.~(\ref{2.4}), we have writing A for the set of points visited by $\pi$,
\begin{equation}\label{2.44}
\begin{array}{lcl}
\limsup\limits_N \;P[\cU^c_{K,uN^d}] & \le & \overline{\lim\limits_N} \;\dsl_{F \in \cL_2} \;\;\dsl_{\pi} \;P[X_{[0,u_0 N^d]} \supseteq A]
\\[1ex]
& \stackrel{(\ref{2.4})}{\le} &  \overline{\lim\limits_N} \;\dsl_{F \in \cL_2}  \;\;\dsl_{\pi} \;\exp \{- \lambda_0 \,|A|\}
\\[1ex]
& \stackrel{(\ref{2.36})}{\le} & \overline{\lim\limits_N} \;  \dsl_{F \in \cL_2} \;8N^27^{[K\log N]-1 }\,e^{-\lambda_0 [K\log N]}
\\
\\[-1ex]
& \le &  \overline{\lim\limits_N} \; c\,N^d \,(7 e^{-\lambda_0})^{[K \log N]},
\end{array}
\end{equation}

\n
where the sum over $\pi$ pertains to the collection of $\star$-nearest neighbor self-avoiding paths with values in F with $[K\log N]$ steps. With (\ref{2.42}) we can thus choose $c_0$ via:

\begin{equation}\label{2.45}
c_0 =  8d \Big(\log \,\mbox{\f $\dis\frac{1}{\mu}$}\Big)^{-1}\,,
\end{equation}
and find
\begin{equation}\label{2.46}
\lim\limits_N \;P[\cU^c_{c_0,uN^d}] = 0, \;\;\mbox{for $u < u_0$}\,,
\end{equation}

\n
from which (\ref{2.40}) follows. We now turn to the proof of (\ref{2.41}). Observe that for $u > 0$, $\ell \ge 1$ and large $N$, one has:
\begin{equation}\label{2.47}
\begin{array}{l}
P\big[\cC^c_{c_0,0,uN^d}\big] \le P\big[X_{[0,uN^d]} \cap B_\infty(0,\ell) \not= \emptyset] \;+
\\[1.5ex]
P\big[X_{[0,uN^d]} \cap B_\infty(0,\ell) = \emptyset, \;\mbox{and} \;\cC^c_{c_0,0,uN^d}\big] \le
\\[1.5ex]
c\,\ell^d\,P\big[0 \in  X_{[0,uN^d]}\big] + \dsl_{F \in \cL_2,0 \in F} \;\dsl_{ \pi} \;P\big[X_{[0,uN^d]} \supseteq A\big]\,,
\end{array}
\end{equation}

\medskip\n
where we have used translation invariance in the last inequality, and the sum over $\pi$ runs over $\star$-nearest neighbor self-avoiding paths with values in $F \cap(B_\infty(0,[ c_0 \log N]) \backslash B_\infty(0,\ell))$, which disconnect $0$ from $F \cap S(0,[c_0 \log N])$, and start on the positive half of the coordinate axis entering the definition of $F$ with smallest label $i \in \{1,\dots,d\}$. As above $A$ stands for the set of points visited by $\pi$. Summing over the different values $k \in \big[\ell +1,[c_0 \log N]\big]$ of the coordinate of the starting point of $\pi$, we see that for small $u$ and sufficiently large $N$, 
\begin{equation}\label{2.48}
\begin{array}{l}
\dsl_{F \in \cL_2, 0 \in F} \;\;\dsl_{\pi} \;P\big[X_{[0,uN^d]} \supseteq A\big] \le c \;\dsl_{k \ge \ell} \;\;\dsl_{m \ge k} \; 7^m \,e^{-\lambda_0 m} =
\\
\\[-1ex]
c\,\dsl_{k \ge \ell} \;(7 e^{-\lambda_0})^{k} (1 - 7 e^{-\lambda_0})^{-1} = c (7 e^{-\lambda_0})^\ell (1 - 7 e^{-\lambda_0})^{-2} \,.
\end{array}
\end{equation}

\medskip\n
Thus coming back to (\ref{2.47}), we see with (\ref{1.13}) that for $u > 0$, $\ell \ge 1$,
\begin{equation}\label{2.49}
\limsup\limits_N \;P[\cC^c_{c_0,0,uN^d}] \le c(1 - e^{-c u}) \,\ell^d + c (7 e^{-\lambda_0})^\ell ( 1 - 7 e^{-\lambda_0})^{-2}\,.
\end{equation}

\n
Letting $u$ tend to $0$ and then $\ell$ to infinity we obtain (\ref{2.41}).
\end{proof}

For $0 < \beta < 1$, and $t \ge 0$, we now introduce the events, cf.~(\ref{0.3}), (\ref{2.37})
\begin{equation}\label{2.50}
\cG_{\beta,t} = \cU_{c_0,t} \cap \cV_{c_0,\beta,t}, \;\mbox{(non-decreasing in $\beta$)}\,.
\end{equation}

\n
The above events encode properties, which enable to single out a giant component. More precisely with the notation of Corollary \ref{cor2.5} we have:

\begin{corollary}\label{cor2.6} $(d \ge d_0, 0 < \beta < 1)$

\medskip
Assume $N \ge 2$, large enough so that $E$ has $|\cdot|_\infty$-diameter bigger than $c_0 \log N$. For $t \ge 0$, on the event $\cG_{\beta,t}$,
\begin{align}
&\mbox{there is a unique connected component in $X^c_{[0,t]}$,  denoted by $O$}, \label{2.51}
\\[-0.5ex]
&\mbox{which contains connected sets $A \in \cA_2$ with $|\cdot|_\infty$-diameter $L_0 \stackrel{\rm def}{=} [c_0 \log N]$}\nonumber
\\[0.5ex] 
&\mbox{(in particular a segment of length $L_0$)}\,,\nonumber
\\[1ex]
&\mbox{for any $F \in \cL_1$, $F \cap O$ contains a segment of length $L_0$}, \label{2.52}
\\[1ex]
&\mbox{the $N^\beta$-neighborhood of $O$ coincides with $E$}\,. \label{2.53}
\end{align}

\medskip\n
Moreover for any $x \in E$,
\begin{equation}\label{2.54}
\mbox{on the event $\cG_{\beta,t} \cap \cC_{c_0,x,t}$, $x$ belongs to $O$}\,.
\end{equation}
Finally one has
\begin{equation}\label{2.55}
\lim\limits_N \;P[\cG_{\beta,uN^d}] = 1, \;\;\mbox{for small $u > 0$}\,.
\end{equation}
\end{corollary}

\begin{proof}
We begin with the proof of (\ref{2.51}) - (\ref{2.53}). With (\ref{0.3}), we see that on $\cG_{\beta,t}$, 
\begin{equation}\label{2.56}
\mbox{any $F \in \cL_1$ contains a segment of length $L_0$ included in $X^c_{[0,t]}$}\,.
\end{equation}

\n
In particular given some $\wt{F} \in \cL_2$, the above applies to all $F \in \cL_1$, with $F \subseteq \wt{F}$. With (\ref{2.37}) any two segments of length $L_0$ in $\wt{F} \backslash X_{[0,t]}$ belong to the same connected component of $F \backslash X_{[0,t]}$ (and hence of $X^c_{[0,t]}$). Now if $\wt{F}$, $\wt{F}_2 \in \cL_2$,
\begin{align}
&\mbox{when $\wt{F}_1 \cap \wt{F}_2 \in \cL_1$, all segments of length $L_0$ in $(\wt{F}_1 \cup \wt{F}_2) \backslash X_{[0,t]}$} \label{2.57}
\\
&\mbox{are in the same connected component of $X^c_{[0,t]}$}\,. \nonumber
\end{align}

\n
Then consider $y \in E$. We can find a nearest neighbor path $(y_i)_{0\le i \le m}$ with $y_0 = 0$, $y_m = y$. Consider $\wt{F} \ni 0$, with $\wt{F} \in \cL_2$, we can construct a sequence $\wt{F}_i \in \cL_2$, $0 \le i \le m$, such that
\begin{align}
&\mbox{$\wt{F}_0 = \wt{F}$, $y_i \in \wt{F}_i$, for $0 \le i \le m$, and either $\wt{F}_{i-1} = \wt{F}_i$ or} \label{2.58}
\\
&\mbox{$\wt{F}_{i-1} \cap \wt{F}_i \in \cL_1$, for $1 \le i \le m$}\,, \nonumber
\end{align}

\n
as we now explain. If $y_1 \in \wt{F}_0 ( = \wt{F})$, we set $\wt{F}_1 = \wt{F}_0$. Otherwise if $y_1 \notin \wt{F}_0$, we choose some canonical vector entering the definition of $\wt{F}_0$ and the canonical vector colinear to $y_1 - y_0$, and define $\wt{F}_1$ as passing through $y_0$ and generated by these two vectors. Clearly $y_1 \in \wt{F}_1$, and $\wt{F}_1 \cap \wt{F}_0 \in \cL_1$. We then continue the construction by induction. 

\medskip
With a similar argument we also see that when $\wt{F},\wt{F}^\prime \in \cL_2$ have a common point $y$ in $E$, we can define $\wt{F}_i \in \cL_2$, $0 \le i \le 2$, such that
\begin{align}
&\mbox{$\wt{F}_0 = \wt{F}$, $\wt{F}_2 = \wt{F}^\prime$, with $y \in \wt{F}_i$, $0 \le i \le 2$, and either $\wt{F}_i = \wt{F}_{i-1}$} \label{2.59}
\\
&\mbox{or $\wt{F}_{i} \cap \wt{F}_{i-1} \in \cL_1$, for $i=1,2$}\,. \nonumber
\end{align}

\n
Combining (\ref{2.56}) - (\ref{2.59}), we see that on $\cG_{\beta,t}$ all segments of length $L_0$ in $X^c_{[0,t]}$ belong to the same connected component of $X^c_{[0,t]}$. With (\ref{2.56}) and the definitions (\ref{0.3}), (\ref{2.37}) of the events entering the definition of $\cG_{\beta,t}$, (\ref{2.51}) - (\ref{2.53}) readily follow. The claim (\ref{2.54}) is a direct consequence of (\ref{2.51}) and (\ref{2.38}). As for (\ref{2.55}) it directly follows from (\ref{1.26}) and (\ref{2.40}).
\end{proof}

In the sequel, on the event $\cG_{\beta,t}$ of (\ref{2.50}), we will refer to the above uniquely defined connected component $O$, as the giant component.

 \section{Excursions to small boxes in a large torus}
 \setcounter{equation}{0}
 
The results of this section are preparatory for the next section, but also of independent interest. We investigate excursions of the random walk to small boxes in the large torus $E = (\IZ / N  \IZ)^d$, with $d \ge 3$. We consider two points $x_1,x_2$ in $E$ with mutual $|\cdot|_\infty$-distance at least $2r + 3$, as well as closed $|\cdot |_\infty$-balls $C(x_i)$, $i = 1,2$, with respective centers $x_i$ and radius $L \le r/10$. We are interested in suitably centered excursions of the walk from the time it first hits $C(x_1) \cup C(x_2)$ up to the last visit to $C(x_1) \cup C(x_2)$ before leaving the closed $r$-neighborhood of $\{x_1,x_2\}$, when the walk is conditioned to leave this $r$-neighborhood at some point $w$ and start at a point $u$ outside this $r$-neighborhood. Of course $w$ determines whether the excursion lies in the neighborhood of $x_1$ or $x_2$, and we center the excursion around $0$ by subtracting the relevant $x_i$ (depending on $w$). As a limit model we consider the excursions of simple random walk on $\IZ^d$ starting with the normalized harmonic measure viewed from infinity of a closed $|\cdot|_\infty$-ball $C$ of radius $L$ and center the origin up to the last visit of $C$. Our main thrust is to derive quantitative controls on the total variation norm between the centered excursions described above and the limit model just explained. Our main result appears in Theorem \ref{theo3.1}. Some of our calculations are similar in spirit to \cite{DembPereRoseZeit06}, see in particular Lemma \ref{lem2.3}. However apart from working in dimension $d \ge 3$, in place of $d=2$, one difference of the results presented here is that they pin-point a limit model for the centered excursions.

\medskip
We now introduce some notation. Throughout this section we assume that $d \ge 3$. We consider positive integers $N,L,r$ such that:
\begin{equation}\label{3.1}
L \ge 1, \;r \ge 10L, \;N \ge 4r + 6\,.
\end{equation}

\n
We define for $x \in E$, (see the beginning of Section 1 for the notation): 
\begin{equation}\label{3.2}
C(x) = B(x,L), \;\wt{C}(x) = B(x,r) \,,
\end{equation}
as well as the subsets of $\IZ^d$,
\begin{equation}\label{3.3}
C  = B(0,L), \;\wt{C}  = B(0,r) \,,
\end{equation}

\n
and tacitly identify $C(0)$ with $C$ and $\wt{C}(0)$ with $\wt{C}$. We then consider two points in $E$.
\begin{equation}\label{3.4}
x_1,x_2 \in E, \;\mbox{with} \;|x_1-x_2|_\infty \ge 2r + 3\,,
\end{equation}

\n
so that $\partial \wt{C}(x_1) \cap \partial \wt{C}(x_2) = \phi$. We then introduce the successive return times to $C(x_1) \cup C(x_2)$ and departures from $\wt{C}(x_1) \cup \wt{C}(x_2)$, cf.~(\ref{1.5}), which we denote with $R_k,D_k, k\ge 1$. In this section we will only need $R_1,D_1$. We also introduce the times of last visits to $C(x_1) \cup C(x_2)$ after $R_k$ and prior to $D_k$:
\begin{equation}\label{3.5}
L_k = \sup\{n \ge R_k, \;X_n \in C(x_1) \cup C(x_2), \;n < D_k\}, \;k \ge 1
\end{equation}

\n
(and for the sake of completeness $L_k$ is defined as $-1$, when the above set is empty, an event which is $P$-negligible). In this section we only consider $L_1$. To describe the centered excursions that interest us, we introduce the canonical space
\begin{align}
\cW = &\; \mbox{the space of finite nearest neighbor $\IZ^d$-valued paths $\ov{w} = (\ov{w}_k)_{0 \le k \le T}$, with} \label{3.6} \\
&\;|\ov{w}_0|_\infty = |\ov{w}_T|_\infty = L\,,\nonumber
\end{align}

\n
denote with $Y_\point$ all the canonical process on $\cW$, and endow the countable space $\cW$ with the $\sigma$-algebra $\cA$ consisting of all subsets of $\cW$.  For $u \notin \wt{C}(x_1) \cup \wt{C}(x_2)$ and $w \in \partial \wt{C}(x_i)$, with $i = 1$ or $2$, we define
\begin{equation}\label{3.7}
\mbox{$Q_{u,w} =$ the law on $\cW$ of $(X_{R_1 + k} - x_i)_{0 \le k \le L_1 - R_1}$ under $P_u[\cdot \,|X_{D_1} = w]$}\,,
\end{equation}

\n
where it should be observed that the conditioning event $\{X_{D_1} = w\}$ has positive probability under $P_u$, and that $P_u [\cdot \,|X_{D_1} = w]$-a.s., $0 < R_1 < L_1 < \infty$ and $X_m \in \wt{C}(x_i)$ for $R_1 \le m \le L_1$, with $i$ as above (\ref{3.7}). So after identification of $\wt{C}$ with $\wt{C}(0)$, (\ref{3.7}) is a meaningful definition.

\medskip
We now turn to the construction of the limit model for these centered excursions. We first introduce the harmonic measure of $C$ viewed from infinity and its mass, the capacity of $C$, cf.~Chapter 2 \S 2 of \cite{Lawl91}:
\begin{align}
e_C(z) & = P_z [\wt{H}_C = \infty], \;\mbox{if $z \in C$, (see (\ref{1.4}) for the notation)}, \label{3.8}
\\
& = 0, \;\mbox{if $z \notin C$}\,, \nonumber
\\[1ex]
{\rm cap}(C) & = e_C(\IZ^d), \;\mbox{and}\label{3.9}
\\[1ex]
\mu_C(z) & = e_C(z) / {\rm cap}(C)\,, \label{3.10}
\end{align}

\n
which is the initial distribution of the limit law. We also define the time of last visit to $C$:
\begin{equation}\label{3.11}
L_C = \sup\{ n \ge 0, \;X_n \in C\}\,,
\end{equation}

\medskip\n
with a similar convention as below (\ref{3.5}) when the above set is empty, and introduce
\begin{equation}\label{3.12}
\mbox{$Q =$ the law on $\cW$ of $(X_k)_{0 \le k \le L_C}$ under $P^{\IZ^d}_{\mu_C}$}\,,
\end{equation}

\n
where $P^{\IZ^d}_{\mu_C}$ stands for the law of simple random walk on $\IZ^d$ with initial distribution $\mu_C$. Note that for any $\ov{w} = (w_k)_{0 \le k \le T}$ in $\cW$,
\begin{equation}\label{3.13}
\begin{split}
Q(Y = \ov{w}) & = E^{\IZ^d}_{\mu_C} [X_k = \ov{w}_k, 0 \le k \le T, \;\mbox{and} \;\wt{H}_C \circ \theta_T = \infty]
\\[1ex]
& = {\rm cap} (C)^{-1} \,e_C (\ov{w}_0) \,P_{\ov{w}_0} [X_k = \ov{w}_k, 0 \le k \le T] \,e_C(\ov{w}_T) \,,
\end{split}
\end{equation}

\n
as a result of the simple Markov property and (\ref{3.8}) - (\ref{3.10}). We are now ready to state the main result of this section.

\begin{theorem}\label{theo3.1} $(d \ge 3)$

\medskip
Assume that {\rm (\ref{3.1}), (\ref{3.4})} hold, and $u \notin \bigcup_{i = 1,2} \wt{C}(x_i)$, $w \in \bigcup_{i=1,2} \partial \wt{C}(x_i)$, then one has:
\begin{equation}\label{3.14}
\|Q_{u,w} - Q\|_{TV} \le c\,\mbox{\f $\dis\frac{L^2}{r}$}\,,
\end{equation}

\n
where for $\nu$ signed measure on $\cW$, $\|\nu\|_{TV} = \sum_{\ov{w} \in \cW} |\nu(\ov{w})|$ denotes  the total variation of $\nu$.
\end{theorem}

\begin{remark}\label{rem3.2} \rm
It will be clear from the proof that the same result holds for collections $x_i$, $1 \le i \le M$, with $|x_i - x_j|_\infty \ge 2r + 3$, whenever $i \not= j$, $u \notin \bigcup_{1 \le i \le M} \wt{C}(x_i)$, $w \in \bigcup_{1 \le i \le M} \partial \wt{C}(x_i)$, (the $\partial \wt{C}(x_i)$, $1 \le i \le M$, are pairwise disjoint due to the above requirement), with $L,r$ as in (\ref{3.1}) and $N \ge M(2r+3)$. As will be clear from the proof below, the constant corresponding to (\ref{3.14}) does not depend on $M$. For simplicity of notation we however restrict to the case $M = 2$. \hfill $\square$
\end{remark}

\begin{proof}
We assume $w \in \partial \wt{C}(x_1)$ and consider $u \notin \wt{C}(x_1) \cup \wt{C}(x_2)$. The case where $w \in \partial \wt{C}(x_2)$ is treated analogously. Note that
\begin{align}
&Q_{u,w} (\wt{\cW}) = 1, \;\mbox{where} \label{3.15}
\\[1ex]
&\wt{\cW} \stackrel{\rm def}{=} \{\ov{w} = (\ov{w}_k)_{0 \le k \le T} \in \cW; \;\;\ov{w}_k \in \wt{C}, \;\mbox{for } \;0 \le k \le T\}\,,\label{3.16}
\end{align}
and that for $\ov{w} \in \wt{\cW}$:
\begin{equation}\label{3.17}
Q_{u,w} (Y = \ov{w}) = A (\ov{w}) \Big/ \Big(\dsl_{\ov{w}^\prime \in \wt{\cW}} A(\ov{w}^\prime)\Big)\,,
\end{equation}
with the notation
\begin{align}
A(\ov{w}) = P_u \big[ &X_{R_1 + k} = x_1 + \ov{w}_k, \,0 \le k \le T, \,X_{n} \notin C(x_1), \;\mbox{for} \label{3.18}
\\
&R_1 + T < n < D_1, \;X_{D_1} = w\big]\,.\nonumber
\end{align}

\n
In what follows when $U$ is a subset of $E$ (resp. $\IZ^d)$, we write $g_{E,U}(\cdot,\cdot)$ (resp. $g_{\IZ^d,U}(\cdot,\cdot)$) to denote the Green function of the walk killed outside $U$, so that
\begin{equation}\label{3.19}
g_{E,U} (x,y) = \dsl_{k \ge 0} \,P_x[X_k = y, k < T_U], \;x,y \in E\,,
\end{equation}

\n
with a similar formula for $g_{\IZ^d,U}$ where $P^{\IZ^d}_x$replaces $P_x$, and $x,y \in \IZ^d$. We simply write $g_{\IZ^d}(\cdot,\cdot)$ when $U = \IZ^d$.

\medskip
Summing over the values of the time of last visit to $(\bigcup_{i=1,2}  \wt{C}(x_i))^c$ before $D_1$, we see that for $\ov{w} \in \wt{W}$:
\begin{equation}\label{3.20}
A(\ov{w}) = \dsl_{v,v^\prime} \;g_{E,\big(\bigcup\limits_{i = 1,2} C(x_i)\big)^c} (u,v) \;\mbox{\f $\dis\frac{1}{2d}$} \;B_{v^\prime}(\ov{w})\,,
\end{equation}

\n
where the above sum runs over $v \sim v^\prime$ with $v \in \partial \wt{C}(x_1)$, $v^\prime \in \wt{C}(x_1)$, with the notation:
\begin{equation}\label{3.21}
\begin{split}
B_{v^\prime}(\ov{w}) & = P_{v^\prime}\big[R_1 + T < T_{\wt{C}(x_1)}, \;X_{R_1+k} = \ov{w}_k + x_1, \;0 \le k \le T\,,
\\
& \qquad \quad X_k \notin C(x_1), \;\mbox{for} \;R_1 + T < k < D_1, \;X_{D_1} = w\big]
\\[1ex]
&=  P^{\IZ^d}_{\wh{v}^\prime} \big[H_C + T < T_{\wt{C}}, \;X_{H_C + k} = \ov{w}_k, \,0 \le k \le T\,, 
\\
& \qquad \quad \ X_k \notin C, \;\mbox{for} \;H_C + T<  k < T_{\wt{C}}, \;X_{T_{\wt{C}}} = \wh{w}\big]\,,
\end{split}
\end{equation}

\n
with the notation $\wh{z} = z - x_1$, using translation invariance and the identification of $\wt{C}(0)$ with $\wt{C}$. Summing over the values of the time of last visit to $\wt{C} \backslash C$ prior to $H_C$, we see that for $\ov{w} \in \wt{\cW}$, $v^\prime$ as above, we find
\begin{equation}\label{3.22}
\begin{split}
B_{v^\prime}(\ov{w})  = &\dsl_{y^\prime} \,g_{\IZ^d, \wt{C} \backslash C} (\wh{v}\,^\prime,y^\prime) \;\dis\frac{1}{2d} \;P_{\ov{w}_0} \big[ X_k = \ov{w}_k, \,0 \le k \le T, \,X_k \notin C ,
\\[-1ex]
&\hspace{4.7cm}  \mbox{for} \; T< k < T_{\wt{C}}, \,X_{T_{\wt{C}}} = \wh{w}\big]
\\[1ex]
= &\dsl_{y^\prime,z^\prime,w^\prime} \,g_{\IZ^d, \wt{C} \backslash C} (\wh{v}\,^\prime,y^\prime) \;\dis\frac{1}{2d} \;P_{\ov{w}_0} \big[ X_k = \ov{w}_k, \,0 \le k \le T\big]  \dis\frac{1}{2d} \;g_{\IZ^d, \wt{C} \backslash C} (z^\prime,w^\prime)\,  \dis\frac{1}{2d}
\end{split}
\end{equation}

\n
where $y^\prime,z^\prime$ run over the respective neighbors in $C^c$ of $\ov{w}_0$ and $\ov{w}_T$, whereas $w^\prime$ runs over the neighbors in $\wt{C}$ of $\wh{w}$, and we have used simple Markov property at times $T+1$, and $T$, and summed over the values of the time of last visit of $\wt{C} \backslash C$ prior to the exit of $\wt{C}$ in $\wh{w}$, to obtain the last expression. The next lemma contains a crucial decoupling estimate.

\begin{lemma}\label{lem3.3} $(d \ge 3$, $L \ge 1$, $10L \le r)$

\medskip
For $a \in \wt{C} \cap \partial (\wt{C}^c), \;b \in \partial C$,
\begin{equation}\label{3.23}
g_{\IZ^d,\wt{C} \backslash C} (a,b) = P^{\IZ^d}_b [T_{\wt{C}} < H_C] \;g_{\IZ^d,\wt{C}}(a,0) (1 + \psi_{a,b})\,,
\end{equation}
where $\psi_{a,b}$ is defined by this equality and
\begin{equation}\label{3.24}
|\psi_{a,b}| \le c_6 \;\dis\frac{L^2}{r}\;.
\end{equation}
\end{lemma}

\begin{proof}
For simplicity we write $g_U(\cdot,\cdot)$ and $g(\cdot,\cdot)$ in place of $g_{\IZ^d,U}(\cdot,\cdot)$ and $g_{\IZ^d}(\cdot,\cdot)$. Using the strong Markov property at time $H_C$, when $H_C < T_{\wt{C}}$, and the symmetry of the killed Green functions, one has
\begin{equation}\label{3.25}
g_{\wt{C}}(a,b) = g_{\wt{C} \backslash C}(a,b) + E^{\IZ^d}_b \big[g_{\wt{C}}(a,X_{H_C}), \,H_C < T_{\wt{C}}]\,.
\end{equation}
Therefore we find:
\begin{equation}\label{3.26}
\begin{split}
g_{\wt{C} \backslash C}(a,b)  = &\;E^{\IZ^d}_b \big[\big(g_{\wt{C}}(a,b) - g_{\wt{C}}(a,X_{H_C})\big), \;H_C < T_{\wt{C}}\big] \;+
\\
& \;g_{\wt{C}}(a,b) \;P^{\IZ^d}_b [H_C > T_{\wt{C}}]
\\[1ex]
= &\; g_{\wt{C}}(a,0) \,P^{\IZ^d}_b [H_C > T_{\wt{C}}] + \big(g_{\wt{C}}(a,b) - g_{\wt{C}} (a,0)\big)\,P^{\IZ^d}_b [H_C > T_{\wt{C}}]\;+
\\
&\;E^{\IZ^d}_b \big[\big(g_{\wt{C}}(a,b) - g_{\wt{C}}(a,X_{H_C})\big), \;H_C < T_{\wt{C}}\big]\,.
\end{split}
\end{equation}

\n
Note that $g_{\wt{C}}(a,\cdot)$ is a non-negative harmonic function on $\wt{C} \backslash \{a\}$. With the Harnack inequality, cf.~Theorem 1.7.2 of \cite{Lawl91}, p.~42,  and a standard covering argument (due to the fact that the quoted theorem refers to Euclidean balls), we find
\begin{equation}\label{3.27}
\sup\limits_{|x|_\infty \le \frac{r}{2}} \;g_{\wt{C}}(a,x) \le c\,g_{\wt{C}}(a,0) \le c\,g(a,0)\,.
\end{equation}

\n
Moreover with the gradient estimates in (a) of Theorem 1.7.1 of \cite{Lawl91}, p.~42, we see that:
\begin{equation}\label{3.28}
\sup\limits_{|x|_\infty \le L,|e| \le 1} \;|g_{\wt{C}}(a,x + e) - g_{\wt{C}}(a,x)| \le \dis\frac{c}{r} \; \sup\limits_{|x|_\infty \le \frac{r}{2}} \;g_{\wt{C}}(a,x)\,.
\end{equation}
Combining (\ref{3.27}), (\ref{3.28}), we see that for all $a \in \wt{C} \cap \partial \wt{C}^c$, 
\begin{equation}\label{3.29}
\sup\limits_{f \in C \cup \partial C} |g_{\wt{C}}(a,f) - g_{\wt{C}}(a,0)| \le c \;\dis\frac{L}{r} \;g_{\wt{C}}(a,0)\,.
\end{equation}

\n
Inserting this inequality in (\ref{3.26}) we see that
\begin{equation}\label{3.30}
\begin{array}{l}
g_{\wt{C} \backslash C} (a,b) - g_{\wt{C}} (a,0) \,P^{\IZ^d}_b [H_C > T_{\wt{C}}]\,\stackrel{\rm def}{=} R, \;\mbox{with}
\\[1ex]
\begin{split}
|R| & \le P^{\IZ^d}_b [H_C > T_{\wt{C}}] \;c\;\dis\frac{L}{r} \;g_{\wt{C}}(a,0) + c\;\dis\frac{L}{r} \;g_{\wt{C}}(a,0)
\\
& \le c\;\dis\frac{L^2}{r} \;g_{\wt{C}}(a,0)\;P^{\IZ^d}_b [H_C > T_{\wt{C}}]\,,
\end{split}
\end{array}
\end{equation}
where in the last step we have used the lower bound
\begin{align*}
P^{\IZ^d}_b [H_C > T_{\wt{C}}] & \ge P^{\IZ^d}_b [H_C > T_{2C}] \times \inf\limits_{x \in (2C)^c} P^{\IZ^d}_x [H_C = \infty]
\\[-0.5ex]
& \ge c\,P^{\IZ^d}_b [H_C > T_{2C}]  \ge c\,P^{\IZ}_1 [H_0 > H_{L+1}] \ge \dis\frac{c}{L}\;.
\end{align*}
Our claim (\ref{3.23}), (\ref{3.24}) now follows.
\end{proof}

We now continue the proof of Theorem \ref{theo3.1}. Note that with the strong Markov property applied at time $H_C$, and standard estimates on the Green function, cf.~\cite{Lawl91}, p.~31, 
\begin{equation}\label{3.31}
\begin{split}
P^{\IZ^d}_z [H_C < \infty] & \le \dsl_{x \in C} \,g_{\IZ^d} (z,x) / \inf\limits_{y \in C} \;\dsl_{x \in C} \,g_{\IZ^d} (y,x)
\\[1ex]
& \le c \,\Big(\dis\frac{L}{r}\Big)^{d-2}, \;\mbox{for} \;z \in \wt{C}^c\,.
\end{split}
\end{equation}

\n
Also by similar estimates as above, and using if necessary the invariance principle to let the path move away, we see with (\ref{3.1}) that
\begin{equation}\label{3.32}
\sup\limits_{z \in \wt{C}^c} \;P^{\IZ^d}_z [H_C < \infty] \le c^\prime < 1\,.
\end{equation}
Hence for $b \in \partial C$, using the strong Markov property at time $T_{\wt{C}}$, we find that
\begin{equation}\label{3.32a}
\begin{split}
0 & \le \;P^{\IZ^d}_b [H_C > T_{\wt{C}}]  - P^{\IZ^d}_b [H_C = \infty]  = P^{\IZ^d}_b [H_C \circ \theta_{T_{\wt{C}}} < \infty, H_C > T_{\wt{C}}]
\\
& \le P^{\IZ^d}_b [H_C > T_{\wt{C}}] \;c^\prime \Big(1 \wedge c  \Big(\dis\frac{L}{r}\Big)^{d-2}\Big), \;\;\mbox{(with $c^\prime < 1$)}\,.
\end{split}
\end{equation}
It thus follows that for $b \in \partial C$,
\begin{equation}\label{3.33}
P^{\IZ^d}_b [H_C > T_{\wt{C}}] = P^{\IZ^d}_b [H_C = \infty] (1 + \e_b),\;\mbox{with $0 \le \e_b \le c \,(\frac{L}{r})^{d-2} \le c_7 \;\frac{L^2}{r}$}\,.
\end{equation}

\n
We now assume for the time being, cf.~(\ref{3.24}), (\ref{3.33}), that
\begin{equation}\label{3.34}
(c_6 + c_7) \;\mbox{\f $\dis\frac{L^2}{r}$} \le \fr \;.
\end{equation}

\n
The case when (\ref{3.34}) does not hold will be straightforwardly handled at the end of the proof. We then define for $a \in \wt{C} \cap \partial (\wt{C}^c)$, $b \in \partial C$, with the notations of (\ref{3.23}), (\ref{3.33}):
\begin{equation}\label{3.35}
\begin{array}{l}
e^{\Gamma_{a,b}} = (1 + \psi_{a,b}) ( 1 + \e_b), \;\mbox{so that} 
\\
|\Gamma_{a,b}| \stackrel{(\ref{3.34})}{\le} c\;\dis\frac{L^2}{r}, \;\mbox{and} \; g_{\IZ^d, \wt{C} \backslash C} (a,b) = P^{\IZ^d}_b [H_C = \infty] \,g_{\IZ^d,\wt{C}}(a,0) \,e^{\Gamma_{a,b}}\,.
\end{array}
\end{equation}

\n
Coming back to (\ref{3.20}), (\ref{3.22}), we see with (\ref{3.23}), (\ref{3.33}) that for $\ov{w} \in  \wt{\cW}$:
\begin{equation}\label{3.36}
\begin{split}
A(\ov{w}) = & \dsl_{y^\prime,z^\prime} \;\Big(\dis\frac{1}{2d}\Big)^2 \;P^{\IZ^d}_{y^\prime} \,[H_C = \infty] \;P^{\IZ^d}_{\ov{w}_0} \,[X_k = \ov{w}_k, 0 \le k \le T] \,P^{\IZ^d}_{z^\prime} [H_C = \infty] \times 
\\[1ex]
&\Big\{\dsl_{v,v^\prime,w^\prime} \;\Big(\dis\frac{1}{2d}\Big)^2\, g_{E,(C(x_1) \cup C(x_2))^c}(u,v) \,g_{\IZ^d,\wt{C}}(v^\prime,0) \,g_{\IZ^d,\wt{C}}(w^\prime,0) \,e^{(\Gamma_{v^\prime,y^\prime} + \Gamma_{w^\prime,z^\prime})}\Big\}\,,
\end{split}
\end{equation}

\n
where in the above sums $y^\prime,z^\prime$ run over $C^c$, with $y^\prime \sim \ov{w}_0$, $z^\prime \sim \ov{w}_T$, $v$ runs over $ \partial \wt{C}(x_1)$, $v^\prime,w^\prime \in \wt{C}$ with $v^\prime \sim \wh{v} ( = v-x_1)$, and $w^\prime \sim \wh{w} = (w - x_1)$. As a result we see that for $\ov{w}_1,\ov{w}_2 \in \wt{\cW}$:
\begin{align}
&\dis\frac{A(\ov{w}_1)}{A(\ov{w}_2)} = \dis\frac{\wt{A}(\ov{w}_1)}{\wt{A}(\ov{w}_2)} \;e^{\Gamma(\ov{w}_1,\ov{w}_2)}, \;\mbox{with} \label{3.37}
\\[1ex] 
&\;\; \,\wt{A}(\ov{w}) = \dsl_{y^\prime,z^\prime}\;\Big(\dis\frac{1}{2d}\Big)^2 \;P^{\IZ^d}_{y^\prime} \,[H_C = \infty] \,P^{\IZ^d}_{\ov{w}_0} [X_k = \ov{w}_k, 0\le k \le T] \,P^{\IZ^d}_{z^\prime} \,[H_C = \infty] \label{3.38}
\\
&\qquad  \;\; \, \stackrel{(\ref{3.8})}{=}  e_C (\ov{w}_0) \,P^{\IZ^d}_{\ov{w}_0} [X_k  = \ov{w}_k, \;0 \le k \le T] \,e_C(\ov{w}_T) \nonumber
\\
&\qquad  \;\; \stackrel{(\ref{3.13})}{=}  {\rm cap}(C) \,Q(Y = \ov{w}), \;\mbox{for $\ov{w} \in \wt{\cW}$, and} \nonumber 
\\[1ex]
&\;\; \;|\Gamma(\ov{w}_1,\ov{w}_2)| \le c \;\dis\frac{L^2}{r}\;. \label{3.39}
\end{align}

\n
Inserting (\ref{3.37}) into (\ref{3.17}), we see that for $\ov{w} \in \wt{\cW}$:
\begin{equation}\label{3.40}
Q_{u,w} [Y = \ov{w}] = \dis\frac{Q(Y = \ov{w})}{\sum\limits_{\ov{w}^\prime \in \wt{\cW}} \,Q(Y = \ov{w}^\prime) \,e^{\Gamma(\ov{w}^\prime,\ov{w})}} = \dis\frac{Q(Y = \ov{w})}{Q(\wt{\cW})} \;e^{G(\ov{w})}\,,
\end{equation}
where $|G(\ov{w})| \le c\;\frac{L^2}{r}$.

\bigskip
Note that with (\ref{3.12}), (\ref{3.16}), and the strong Markov property
\begin{equation}\label{3.41}
Q(\wt{\cW}^c) = P^{\IZ^d}_{\mu_C} [H_C \circ \theta_{T_{\wt{C}}} < \infty] \stackrel{(\ref{3.31}),(\ref{3.32})}{\le} c^\prime \Big(1 \wedge c \Big(\dis\frac{L}{r}\Big)^{d-2}\Big), \;\mbox{with $c^\prime < 1$}\,.
\end{equation}
We thus find that
\begin{equation}\label{3.42}
\begin{split}
\|Q_{u,w} - Q\|_{TV} &= \dsl_{\ov{w} \in \wt{\cW}} \,|Q_{u,w} (Y = \ov{w}) - Q(Y = \ov{w})| + Q(\wt{\cW}^c)
\\
&=  \dsl_{\ov{w} \in \wt{\cW}} \;\dis\frac{Q(Y = \ov{w})}{Q(\wt{\cW})} \;|\exp\{G(\ov{w})\} -1 + Q(\wt{\cW}^c) | + Q(\wt{\cW}^c) \le c\;\dis\frac{L^2}{r}\;,
\end{split}
\end{equation}
using (\ref{3.40}), (\ref{3.41}). As a result we have proved (\ref{3.14}) under (\ref{3.34}). On the other hand when (\ref{3.34}) does not hold, $\frac{L^2}{r} \ge \frac{1}{2}\, (c_6 + c_7)^{-1}$, and
\begin{equation*}
\|Q_{u,w} - Q\|_{TV} \le 2 \le 4(c_6 + c_7) \;\dis\frac{L^2}{r} \,,
\end{equation*}

\n
so that adjusting the constant in (\ref{3.14}) if necessary, we have completed the proof of Theorem \ref{theo3.1}.
\end{proof}

\section{Volume estimate for the giant component}
 \setcounter{equation}{0}
 
The main purpose of this section is to show that the giant component $O$ in the vacant set left by the walk at time $t = u N^d$, (this component is well-defined on the event $\cG_{\beta,t}$, cf.~(\ref{2.51})), typically occupies a non-degenerate fraction of the volume of the torus $E$, when $N$ is large and $u$ chosen small. The statement (\ref{2.54}) provides a local criterion, depending on the configuration of vacant sites left by the walk in a neighborhood of order const $\log N$ of a point $x \in E$, which ensures, when $\cG_{\beta,t}$ occurs, that $x$ belongs to $O$. With (\ref{2.55}) this reduces the problem of proving the non-degeneracy of the volume of $O$ to the question of showing that typically the asymptotic fraction of points $x$ in $E$ that fulfill the local condition $\cC_{c_0,x,uN^d}$ is non-degenerate when $u$ is small. With (\ref{2.4}) this task is further reduced to the control on the variance of this quantity. As it turns out it is simpler to bound the variance of the fraction of points of $E$ that satisfy a modified local condition where the fixed time $t = uN^d$, is replaced by a random time corresponding to the completion of const $u(\log N)^{2(d-2)}$ excursions of the walk to a neighborhood of order $(\log N)^2$ of the point, cf.~(\ref{4.22}). The controls of Section 3 are then instrumental in bounding the variance of this modified ratio, cf.~Proposition \ref{prop4.2}. Our main estimates on averages of suitable local functions are expressed in a general form, (not specifically referring to (\ref{2.54})), and appear in Theorem \ref{theo4.3}, when $d \ge 3$. The applications to the vacant set, the giant component (when $d \ge d_0$, cf.~Corollary \ref{cor2.5}), and the size of the largest ball in the vacant set are given in Corollary \ref{cor4.5}, \ref{cor4.6}, \ref{cor4.8}.
 
 \medskip
We now begin with some additional notation. We consider $d \ge 3$, $L \ge 1$, $r \ge 10L$, $N \ge 10r$, $x \in E$, and recall the definition of $C(x) \subseteq \wt{C}(x)$ in (\ref{3.2}). We consider some function $\phi$, defined on the collection of subsets of $C(0)$:
 \begin{equation}\label{4.1}
 \phi: A \subseteq C(0) \rightarrow \phi(A) \in [0,1] \,.
 \end{equation}
 
 \n
 Typical examples to keep in mind are for instance
 \begin{align}
 \phi_0(A) & \stackrel{\rm def}{=} 1 \{0 \notin A\}, \;\mbox{for $A \subseteq C(0)$, or} \label{4.2}
 \\[1ex]
 \phi_1(A) &\, = 1\{\mbox{for some $F \in \cL_2$, with $0 \in F$, $0$ is connected to $S(0,L)$ in $F \backslash A\}$}\,, \label{4.3}
 \end{align}
 
 \n
where we recall (\ref{2.38}) for the latter example.With $\phi$ as in (\ref{4.1}), we then define for $x \in E$ and $t \ge 0$:
\begin{equation}\label{4.4}
 h(x,t) = \phi \big(\big(X_{[0,t]} \cap C(x)\big) - x\big)\,.
 \end{equation}
 
\medskip \n
Our chief task in this section consists in the derivation of appropriate lower bounds on ratios of the type:
\begin{equation}\label{4.5}
\Gamma_u = \dis\frac{1}{N^d} \;\dsl_{x \in E} \;h(x,uN^d), \;\;\mbox{with $u > 0$}\;.
\end{equation}

\n
For $x \in E$, we introduce in analogy to (\ref{1.9}), (\ref{1.10}),
\begin{equation}\label{4.6}
B(x) = x + B \subseteq \wt{B}(x) = x + \wt{B}, \;
(\mbox{so}\; C(x) \varsubsetneq \wt{C}(x) \varsubsetneq B(x) \varsubsetneq \wt{B}(x))
\end{equation}

\n
as well as the successive returns to $B(x)$ and departures from $\wt{B}(x)$: 
\begin{equation}\label{4.7}
\cR^x_k, \cD^x_k, \,k \ge 1\,.
\end{equation}

\n
We also consider, cf.~(\ref{1.5}), the successive returns to $C(x)$ and departures from $\wt{C}(x)$: 
\begin{equation}\label{4.8}
R^x_k,D^x_k,k \ge 1\,.
\end{equation}

\n
We begin with the following auxiliary result,(note that $r$ does not appear in the right-hand side of the inequalities):
\begin{lemma}\label{lem4.1} $(d \ge 3, L \ge 1$, $r \ge 10L$, $N \ge 10r$)

\medskip
There are constants $c_8 > c_9 > 0$, such that for $u > 0$, $x \in E$:
\begin{align}
P[R^x_{\ell^*(u)} \le uN^d] \le c \,e^{-cu\,L^{d-2}}, &\;\, \mbox{with} \;\,\ell^* (u) = [c_8 \,u\,L^{d-2}]\,,\label{4.9}
\\[1ex]
P[D^x_{\ell_*(u)} \ge uN^d] \le c \,e^{-cu\,L^{d-2}}, &\;\, \mbox{with} \;\,\ell_* (u) = [c_9 \,u\,L^{d-2}]\,.\label{4.10}
\end{align}
\end{lemma}

\begin{proof}
We begin with the proof of (\ref{4.9}). We introduce for $\ell \ge 1$
\begin{equation}\label{4.11}
Z^x_\ell = \dsl_{m \ge 1} \;1 \{\cR^x_\ell \le R^x_m \le \cD_\ell^x\} = \dsl_{m \ge 1} 1\{\cR^x_\ell \le D^x_m \le \cD^x_\ell\}\,.
\end{equation}

\n
With the strong Markov property at times $D^x_m$ and $H_{\wt{C}(x)}$, we see that for $i \ge 0$, $\ell \ge 2$, $P$-a.s.,
\begin{equation}\label{4.12}
\begin{array}{l}
P[Z^x_\ell > i\,|\,\cF_{\cR^x_\ell}] = P_{X_{\cR^x_\ell}} [R^x_{i+1} < T_{\wt{B}(x)}]\le 
\\[1ex]
P_{X_{\cR^x_\ell}} [H_{\wt{C}(x)} < T_{\wt{B}(x)}] \,\big(\sup\limits_{|z - x|_\infty \in \{r,r+1\}} P_z[H_{C(x)} < T_{\wt{B}(x)}]\big)^{i+1} =
\\
\\[-2ex]
P_{X_{\cR^x_\ell}} [H_{\wt{C}(x)} < T_{\wt{B}(x)}] \;\big(\sup\limits_{|z - x|_\infty = r} \;P_z[H_{C(x)} < T_{\wt{B}(x)}])^{i+1} \,.
\end{array}
\end{equation}

\n
Analogously we have for $i \ge 0$
\begin{equation}\label{4.13}
P[Z_1^x > i] \le P[H_C < T_{\wt{B}}] \;\big(\sup\limits_{|z- x|_\infty = r} P_z[H_{C(x)} < T_{\wt{B}(x)}]\big)^i\,.
\end{equation}

\n
Using similar bounds as in (\ref{3.31}), (\ref{3.32}), we find that for $\ell \ge 2$, $i \ge 0$, $P$-a.s.,
\begin{equation*}
P[Z^x_\ell > i \,| \,\cF_{\cR^x_\ell}] \le c^\prime \wedge \Big(c \Big(\dis\frac{r}{N}\Big)^{d-2}\Big) \cdot \Big\{ c^\prime \wedge \Big(c \Big(\dis\frac{L}{r}\Big)^{d-2}\Big) \Big\}^{i+1} \;\mbox{with $c^\prime < 1$}\,.
\end{equation*}
Using the inequality
\begin{equation}\label{4.14}
\begin{array}{l}
P[H_{C(x)} < T_{\wt{B}(x)}] \le 
\\[1ex]
E\Big[\dsl_{k \ge 0} \,1\{X_k \in C(x), k < T_{\wt{B}(x)}\}\Big]\; \big/ \inf\limits_{y \in C(x)} E_y \Big[\dsl_{k \ge 0} \,1\{X_k \in C(x), k < T_{\wt{B}(x)}\}\Big]\,,
\end{array}
\end{equation}

\n
a similar upper bound as in (\ref{1.34}) on the numerator and a lower bound of type $c\,L^2$ on the denominator with the help of the invariance principle, we find that
\begin{equation*}
P[H_{C(x)} < T_{\wt{B}(x)}] \le c\,\Big(\mbox{\f $\dis\frac{L}{N}$}\Big)^{d-2} \,,
\end{equation*}

\n
and it is also straightforward to argue with the invariance principle and similar arguments as for the derivation of (\ref{3.32}) that the above probability is bounded by some $c^{\prime\prime} < 1$. Coming back to (\ref{4.12}), (\ref{4.13}), we thus see that:
\begin{equation}\label{4.15}
\begin{array}{l}
P[Z^x_\ell > i\,| \,\cF_{\cR^x_\ell}] \le p_0 \,p^{i+1}, \;\;P[Z^x_1 > i] \le p_0 \,p^{i+1}, \;\mbox{for $i \ge 0$, $\ell \ge 2$},
\\[1ex]
\mbox{with $p_0 = c_{10} \wedge \Big(c \,\Big(\dis\frac{r}{N}\Big)^{d-2}\Big), \quad p = c_{10} \wedge  \Big(c \,\Big(\dis\frac{L}{r}\Big)^{d-2}\Big)$, with $c_{10} < 1$}\,.
\end{array}
\end{equation}

\n
With an argument of stochastic domination, we thus see that for $\lambda > 0$ with $e^\lambda \,p < 1$, and $\ell \ge 2$,
\begin{equation}\label{4.16}
\begin{array}{l}
E[\exp \{\lambda\,Z^x_\ell\} \,|\,\cF_{\cR^x_\ell}] \le 1 - p_0 \,p + \dsl_{k \ge 1} \,e^{\lambda k} p_0 \,p^k (1-p) = 1 + p_0 \,p\;\dis\frac{(e^\lambda - 1)}{(1-e^\lambda p)} \;,
\\[1ex]
E[\exp \{\lambda\,Z^x_1\}] \le 1 + p_0 \,p \;\dis\frac{(e^\lambda - 1)}{(1- e^\lambda p)} \;.
\end{array}
\end{equation}

\n
As a result we find that with the notation below (\ref{1.12})
\begin{equation}\label{4.17}
\begin{split}
P[Z_1 + \dots + Z_{k^*} \ge n] & \le \exp \{-\lambda n\} \;\Big(1 + p_0 \,p \,\dis\frac{(e^\lambda - 1)}{1 - e^\lambda\,p}\Big)^{k^*} 
\\[1ex]
&  \le \exp \Big\{ - \lambda n + k^* p_0 \,p \,\dis\frac{(e^\lambda - 1)}{1 - e^\lambda p} \Big\}\,.
\end{split}
\end{equation}

\n
Note that $k^* p_0 \,p \le c \,u\,L^{d-2}$, and choosing $\lambda$ so that $e^\lambda c_{10} = \frac{1}{2} \;(1 + c_{10})$, (recall $p \le c_{10} < 1)$, we thus obtain:
\begin{equation}\label{4.18}
\begin{array}{l}
P[R^x_n \le u N^d] \le P[\cR^x_{k^*} \le u N^d] + P [\cR^x_{k^*} > u N^d \ge R^x_n]
\\[1ex]
\qquad \qquad \;\; \stackrel{(\ref{1.11}),(\ref{4.17})}{\le} c\,\exp\{- c u N^{d-2}\} + c\,\exp\{ - \lambda n + cu L^{d-2}\}\,,
\end{array}
\end{equation}

\n
and (\ref{4.9}) follows straightforwardly.

\medskip
We now turn to the proof of (\ref{4.10}). We use a bound from below on $P[H_{C(x)} < T_{\wt{B}(x)}]$ and $P_z[H_{C(x)} < T_{\wt{B}(x)}]$ with a similar right-hand side as in (\ref{4.14}), except for the fact that inf is replaced with sup, and in the case of the second probability $E$ is replaced with $E_z$, see also (\ref{1.57}). Then with standard Green function estimates, see for instance (1.11) of \cite{DembSzni06}, we obtain:
\begin{equation}\label{4.19}
P[Z^x_\ell > 0 \,| \,\cF_{\cR^x_\ell}] \stackrel{P{\rm -a.s.}}{\ge}  c \,\Big(\dis\frac{L}{N}\Big)^{d-2},    \;\mbox{for $\ell \ge 2$, and $P[Z_1^x > 0] \ge c\,\Big(\dis\frac{L}{N}\Big)^{d-2}$}\,.  
\end{equation}

\n
As a result we see that for $\lambda > 0$, $\ell \ge 2$,
\begin{align*}
& E[\exp\{- \lambda Z^x_\ell\} \,|\,\cF_{\cR^x_\ell}] \le 1 - (1-e^{-\lambda}) \,c \,\Big(\dis\frac{L}{N}\Big)^{d-2}, \;\mbox{and}
\\[1ex]
&E[\exp\{-\lambda Z_1^x\}] \le 1 - (1 - e^{-\lambda}) \,c\,\Big(\dis\frac{L}{N}\Big)^{d-2}\,,
\end{align*}
so that for $n \ge 1$, (with the convention that the sum in the probability below vanishes when $k_* \le 1$),
\begin{equation}\label{4.20}
P[Z_1 + \dots + Z_{(k_* - 1)_+} < n] \le \exp\Big\{ \lambda n - (k_* - 1)_+ (1 - e^{-\lambda}) \,c\,\Big(\dis\frac{L}{N}\Big)^{d-2}\Big\}\,,
\end{equation}
where $(k_* - 1)_+ \,c(\frac{L}{N})^{d-2} \ge c \,u\,L^{d-2} - c$, cf.~below (\ref{1.12}). We then see that
\begin{equation}\label{4.21}
\begin{array}{l}
P\big[D^x_n \ge u N^d\big] \le P[\cR^x_{k_*} \ge uN^d] + P[\cR^x_{k_*} < uN^d, \;D^x_n \ge \cR^x_{k_*}]
\\[0.5ex]
\stackrel{(\ref{1.12})}{\le} \;c\;\exp \{- c u\,N^{d-2}\} + P[Z_1 + \dots + Z_{(k_*-1)_+} < n]
\\[0.5ex]
 \stackrel{(\ref{4.20})}{\le} c\,\exp \{ - c u\,N^{d-2}\} + c\,\exp\{\lambda n - (1-e^{-\lambda}) \,c \,u L^{d-2}\}\,.
\end{array}
\end{equation}

\n
Choosing $\lambda$ so that $e^{-\lambda} = \frac{1}{2}$, the claim (\ref{4.10}) follows straightforwardly.
\end{proof}

We now introduce a modification of $\Gamma_u$ in (\ref{4.5}), which is more convenient when bounding its variance. Namely we define with (\ref{4.1}),  (\ref{4.4}), and the notation from (\ref{4.9})
\begin{equation}\label{4.22}
\wt{\Gamma}_u = \dis\frac{1}{N^d} \;\dsl_{x \in E} \;h(x, D^x_{\ell^*(u)}), \;\mbox{for $u > 0$}\,.
\end{equation}

\n
Our main estimate on the variance of $\wt{\Gamma}_u$ comes in the next proposition. In what follows var and cov denote the variance and covariance under $P$.

\begin{proposition}\label{prop4.2} ($d \ge 3$, $L \ge 1$, $N \ge 10r$, $r \ge 10L$, under (\ref{4.1}))
\begin{equation}\label{4.23}
{\rm var}(\wt{\Gamma}_u) \le c((\frac{r}{N})^d + u \,\frac{L^d}{r}), \;\mbox{for $u > 0$}\,.
\end{equation}
\end{proposition}

\begin{proof} When $\ell^*(u) = 0$, with our conventions we see that $\wt{\Gamma}_u = \frac{1}{N^d} \,\sum_{x \in E} \phi(C(x) \cap \{X_0\} - x)$, a non-random quantity as follows from translation invariance. The claim (\ref{4.23}) is then trivially satisfied. We thus assume from now on that $\ell^*(u) \ge 1$.
We then consider an integer $r$ as in (\ref{4.23}), and write
\begin{equation}\label{4.24}
\begin{split}
{\rm var} (\wt{\Gamma}_u) & = \dis\frac{1}{N^{2d}} \;\dsl_{x_1,x_2 \in E} \;{\rm cov} \big(h (x_1,D^{x_1}_{\ell^*(u)}), \;h(x_2,D^{x_2}_{\ell^*(u)})\big)
\\[1ex]
& \le c \,\Big(\dis\frac{r}{N}\Big)^d + \sup\limits_{|x_1-x_2|_\infty \ge 2r+3} \big|{\rm cov}\big(h(x_1,D^{x_1}_{\ell^*(u)}), \;h(x_2,D^{x_2}_{\ell^*(u)})\big)\big|\,.
\end{split}
\end{equation}

\n
We recall the notations $R_k,D_k,k \ge 1$ introduced below (\ref{3.4}), and write for $i = 1,2$:
\begin{equation}\label{4.25}
\begin{split}
n^{x_i}_1 & = \inf\{k \ge 1; \;X_{R_k} \in C(x_i)\}, \;\mbox{and}
\\[1ex]
n^{x_i}_{j+1} & = \inf\{k > n_j^{x_i}; \;X_{R_k} \in C(x_i)\}, \;j \ge 1\,,
\end{split}
\end{equation}

\medskip\n
The relation between $R^{x_i}_k$, $D^{x_i}_k$, $k \ge 1$, for $i = 1,2$, and $R_k, D_k, k \ge 1$, is the following: one has $P$-a.s.,
\begin{equation}\label{4.26}
R^{x_i}_k = R_{n^{x_i}_k}, \;\;D^{x_i}_k = D_{n^{x_i}_k}, \;\;\mbox{for $k \ge 1$, $i = 1,2$}\,.
\end{equation}
We then introduce the constant, cf.~(\ref{4.9}), (\ref{4.10})
\begin{equation}\label{4.27}
c_{11} = \dis\frac{c_8}{c_9} > 1\,.
\end{equation}

\n
We recall the definition (\ref{3.5}) and denote with $e_k(\cdot)$ the $P$-a.s.~well-defined centered excursion
\begin{equation}\label{4.32}
e_k(m) = X_{R_k + m} - x_i, \;0 \le m \le L_k - R_k, \;\mbox{on} \; \{X_{R_k} \in C(x_i)\}, \,i=1,2\,.
\end{equation}

\n
We recall our tacit identification of $\wt{C}(0) \subseteq E$ with $\wt{C}$ in $\IZ^d$, see below (\ref{3.3}), so that $P$-a.s., $e_k(\cdot) \in \wt{\cW} \subseteq \cW$, cf.~(\ref{3.6}), (\ref{3.16}). We also consider the $k$-th excursion to $C(x_i)$, after centering at the origin, which is also $P$-a.s. well-defined:
\begin{equation}\label{4.33}
e^i_k(\cdot) = e_{n_k^{x_i}}(\cdot) \;k \ge 1, \;\;i \in \{1,2\}\,,
\end{equation}
as well as its trace
\begin{equation}\label{4.34}
\cS^i_k = {\rm Im} \;e^i_k\,,
\end{equation}

\n
where for $\ov{w} = (\ov{w}_m)_{0 \le m \le T} \in \cW, \;Im \,\ov{w} = \{w_0,\dots,w_T\} \subseteq \IZ^d$. With the above notation, we see that $P$-a.s.,
\begin{align}
X_{[0,D^{x_i}_{\ell^*(u)}]} \cap C(x_i) - x_i & = (\cS^i_1 \cup \dots \cup \cS^i_{\ell^*(u)} ) \cap C, \;\mbox{and} \label{4.35}
\\[1ex]
h(x_i, D^{x_i}_{\ell^*(u)}) & = G(e^i_1,\dots,e^i_{\ell^*(u)}), \;\mbox{for $i = 1,2$},\label{4.36}
\end{align}

\n
where $G$ is the function from $\cW^{\ell^*(u)}$ into $[0,1]$ defined by, cf.~(\ref{4.1}),
\begin{equation}\label{4.37}
G(\ov{w}_1,\dots,\ov{w}_{\ell^*(u)}) = \phi\big((Im \,\ov{w}_1 \cup \dots \cup Im \,\ov{w}_{\ell^*(u)}) \cap C\big) \,.
\end{equation}

\n
We now consider two $[0,1]$-valued functions $G_1$, $G_2$ on $\cW^{\ell^*(u)}$, (we are especially interested in the case $G_i = G$ or $G_i = 1$), and write
\begin{equation}\label{4.38}
H_i = G_i (e^i_1,\dots,e^i_{\ell^*(u)}), \;i = 1,2 \,.
\end{equation}
We see that for $z \notin \wt{C}(x_1) \cup \wt{C}(x_2)$
\begin{equation}\label{4.39}
E_z[H_1H_2] = \dsl_{\cK} E_z[H_1 H_2, \,A_{(\ov{k}{\hspace{-0.2ex}\,^1},\ov{k}{\hspace{-0.2ex}\,^2})}]\,,
\end{equation}
where $\cK$ denotes the set of ordered pairs of $\ell^*(u)$-uples of integers
\begin{align}
&\mbox{$(\ov{k}{\hspace{-0.2ex}\,^1},\ov{k}{\hspace{-0.2ex}\,^2})$, with $1 \le \ov{k}{\hspace{-0.2ex}\,^i_1} < \dots < \ov{k}{\hspace{-0.2ex}\,^i_{\ell^*(u)}}$, for $i = 1,2$}, \label{4.40}
\\
&\mbox{with all $\ov{k}{\hspace{-0.2ex}\,^i_j}$ distinct, for $1 \le j \le \ell^*(u)$, $i = 1,2$,} \nonumber
\end{align}
and for $(\ov{k}{\hspace{-0.2ex}\,^1},\ov{k}{\hspace{-0.2ex}\,^2}) \in \cK$, we write
\begin{equation}\label{4.41}
A_{(\ov{k}{\hspace{-0.2ex}\,^1},\ov{k}{\hspace{-0.2ex}\,^2})} = \big\{n^{x_i}_m =  \ov{k}{\hspace{-0.2ex}\,^i_m}, \;\mbox{for $1 \le m \le \ell^*(u), \,i = 1,2\big\}$}\,.
\end{equation}
We introduce the $\sigma$-algebra
\begin{equation}\label{4.42}
\mbox{$\cE =$ the $P$-completion of $\sigma(X_{D_k}, k \ge 1)$}\,.
\end{equation}
Note that $n^{x_i}_m$, $i \in \{1,2\}$, $m \ge 1$, are $\cE$-measurable, so that one has 
\begin{equation}\label{4.43}
A_{(\ov{k}{\hspace{-0.2ex}\,^1},\ov{k}{\hspace{-0.2ex}\,^2})} \in \cE, \;\mbox{for any} \;(\ov{k}{\hspace{-0.2ex}\,^1},\ov{k}{\hspace{-0.2ex}\,^2})\in \cK\,.
\end{equation}
Using the strong Markov property at the times $D_m, m \le \ov{k} \stackrel{\rm def}{=} \ov{k}{\hspace{-0.2ex}\,^1_{\ell^*(u)}} \vee \ov{k}{\hspace{-0.2ex}\,^2_{\ell^*(u)}}$, we see that for $(\ov{k}{\hspace{-0.2ex}\,^1}, \ov{k}{\hspace{-0.2ex}\,^2}) \in \cK$, $P$-a.s. on $A_{(\ov{k}{\hspace{-0.2ex}\,^1},\ov{k}{\hspace{-0.2ex}\,^2})}$,
\begin{equation}\label{4.44}
\begin{array}{l}
E_z [ H_1 H_2 \,| \cE] = 
\\
\dis\int \,G_1 (\ov{w}_{\ov{k}{\hspace{-0.2ex}\,^1_1}}, \dots , \ov{w}_{\ov{k}{\hspace{-0.2ex}\,^1_{\ell^*(u)}}}) \, G_2 (\ov{w}_{\ov{k}{\hspace{-0.2ex}\,^2_1}}, \dots , \ov{w}_{\ov{k}{\hspace{-0.2ex}\,^2_{\ell^*(u)}}}) \;\prod\limits^{\ov{k}}_{m=1} \,Q_{X_{D_{m-1}}},_{X_{D_m}} (d \ov{w}_m)\,,
\end{array}
\end{equation}

\n
where we used the notation (\ref{3.7}), and the convention $X_{\hspace{-0.4ex} D_0} = z$, when $m = 1$.

\medskip
We can now find for each $u \notin \wt{C}(x_1) \cup \wt{C}(x_2)$, $w \in \partial \wt{C}(x_1) \cup \partial \wt{C}(x_2)$ a coupling $\wt{Q}_{u,w} (d\ov{w}, d\ov{w}{\,^\prime})$ on $\cW \times \cW$ such that, see (\ref{3.12}),
\begin{align}
&\mbox{under the first (resp. the second) canonical coordinate the image of} \label{4.45}
\\
&\mbox{$\wt{Q}_{u,w}$ is $Q_{u,w}$ (resp. $Q$)}\,,\nonumber
\\[1ex]
&\wt{Q}_{u,w} (\ov{w} \not= \ov{w}{\,^\prime}) = \fr \;\|Q_{u,w} - Q \|_{TV} \stackrel{(\ref{3.14})}{\le} c\;\dis\frac{L^2}{r}\;,\label{4.46}
\end{align}

\n
for the construction of $\wt{Q}_{u,w}$ see for instance Theorem 5.2, p.~19 of \cite{Lind92}. We thus see that for $(\ov{k}{\hspace{-0.2ex}\,^1}, \ov{k}{\hspace{-0.2ex}\,^2}) \in \cK$, $P$-a.s. on $A_{(\ov{k}{\hspace{-0.2ex}\,^1},\ov{k}{\hspace{-0.2ex}\,^2})}$,
\begin{equation}\label{4.47}
\begin{array}{l}
|E_z [ H_1 H_2 \,| \cE]  - \prod\limits^2_{i=1} \;E^{Q^{\otimes \ell^*(u)}}[G_i] \,| =
\\
\Big| \dis\int \,G_1 (\ov{w}_{\ov{k}{\hspace{-0.2ex}\,^1_1}}, \dots , \ov{w}_{\ov{k}{\hspace{-0.2ex}\,^1_{\ell^*(u)}}}) \, G_2 (\ov{w}_{\ov{k}{\hspace{-0.2ex}\,^2_1}}, \dots , \ov{w}_{\ov{k}{\hspace{-0.2ex}\,^2_{\ell^*(u)}}}) \; - 
\\[2ex]
G_1 (\ov{w}^{\,\prime}_{\ov{k}{\hspace{-0.2ex}\,^1_1}}, \dots , \ov{w}^{\,\prime}_{\ov{k}{\hspace{-0.2ex}\,^1_{\ell^*(u)}}}) \, G_2 (\ov{w}^{\,\prime}_{\ov{k}{\hspace{-0.2ex}\,^2_1}}, \dots , \ov{w}^{\,\prime}_{\ov{k}{\hspace{-0.2ex}\,^2_{\ell^*(u)}}}) \;\prod\limits^{\ov{k}}_{m=1} \,\wt{Q}_{X_{D_{m-1}},X_{D_m}} (d \ov{w}_m, d\ov{w}_m^{\,\prime}) \Big| \le
\\[2ex]
2 \ell^*(u) \sup\limits_{m \in \{\ov{k}{\hspace{-0.2ex}\,^i_j};\, i = 1,2, \,1 \le j \le \ell^*(u)\}} \wt{Q}_{X_{D_{m-1}},X_{D_m}} (\ov{w}_m \not= \ov{w}^{\,\prime}_m) \stackrel{(\ref{4.46})}{\le} c\ell^*(u) \,\dis\frac{L^2}{r} \stackrel{(\ref{4.9})}{\le} cu\,\dis\frac{L^d}{r}\,.
\end{array}
\end{equation}

\n
Hence with (\ref{4.39}), (\ref{4.43}), we see that for $z \notin \wt{C}(x_1) \cup \wt{C}(x_2)$:
\begin{equation*}
|E_z[H_1 H_2] - \prod\limits^2_{i=1} \;E^{Q^{\otimes \ell^*(u)}} [G_i]| \le c  u\,\dis\frac{L^d}{r}\;,
\end{equation*}
and hence
\begin{equation}\label{4.48}
|E[H_1 H_2] - \prod\limits^2_{i=1} \;E^{Q^{\otimes \ell^*(u)}}[G_i]| \le c\,\Big(\Big(\dis\frac{r}{N}\Big)^d + u\;\dis\frac{L^d}{r}\Big)
\end{equation}

\n
choosing $G_i = G$, cf.~(\ref{4.37}), or $G_i = 1$, we see with (\ref{4.36}) that the last term in the second line of (\ref{4.24}) is smaller than $c((\frac{r}{N})^d + u \,\frac{L^d}{r})$. With (\ref{4.24}), the claim (\ref{4.23}) follows.
\end{proof}

When the function $\phi$ in (\ref{4.1}) is monotone decreasing, i.e. for $A \subseteq A^\prime \subseteq C(0)$, $\phi(A) \ge \phi(A^\prime)$, then we can easily transfer controls on $\Gamma_\point$ from controls on $\wt{\Gamma}$.

\begin{theorem}\label{theo4.3} $(d \ge 3, L \ge 1, N \ge 100L)$

\medskip
Assume that $\phi$ in(\ref{4.1}) is monotone decreasing, then for $u > 0, s > 0$,
\begin{equation}\label{4.49}
P[\Gamma_u < E[\Gamma_{c_{11} u}] - c\,\exp\{- c\,u \,L^{d-2}\} - s] \le c\;\dis\frac{\sigma^2_{u,L,N}}{s^2} + c\,N^d\exp\{ - c\,u \,L^{d-2}\}\,,
\end{equation}
and
\begin{equation}\label{4.50}
P[\Gamma_u > E[\Gamma_{c^{-1}_{11} u}] + c\,\exp\{- c u \,L^{d-2}\} + s] \le c\;\dis\frac{\sigma^2_{u,L,N}}{s^2} + c\,N^d\exp\{ - c u \,L^{d-2}\}\,,
\end{equation}

\n
where $c_{11} > 1$ is defined in (\ref{4.27}) and
\begin{equation}\label{4.51}
\sigma^2_{u,L,N} \stackrel{\rm def}{=} \inf\Big\{\Big(\dis\frac{r}{N}\Big)^d + u\;\dis\frac{L^d}{r} \;; \;10 L \le r \le \dis\frac{N}{10}\Big\}\,.
\end{equation}
\end{theorem}

\begin{proof}
Choose $r$ as in (\ref{4.51}) and define $\wt{\Gamma}_u$ as in (\ref{4.22}). Since $\phi$ is monotone decreasing, we see that
\begin{align*}
E [\Gamma_u] - E[\wt{\Gamma}_u]  & = \dis\frac{1}{N^d} \;\dsl_{x \in E} \;E[h(x,uN^d) - h(x, D^x_{\ell^*(u)})] 
\\
& \ge - \dis\frac{1}{N^d} \;\dsl_{x \in E} \;P[D^x_{\ell^*(u)} < uN^d] \stackrel{\rm (\ref{4.9})}{\ge} - c\,e^{-c  u\,L^{d-2}}\,,
\end{align*}

\n
and using the fact that $\ell^*(\frac{u}{c_{11}}) = \ell_* (u)$, cf.~(\ref{4.27}), (\ref{4.9}), (\ref{4.10}), we also have
\begin{align*}
E [\wt{\Gamma}_{\frac{u}{c_{11}}}] - E[\Gamma_u]  & = \dis\frac{1}{N^d} \;\dsl_{x \in E} \; E[h(x,D^x_{\ell_*(u)}) - h(x,uN^d)]
\\[1ex]
& \ge  - \dis\frac{1}{N^d} \;\dsl_{x \in E} \; P[D^x_{\ell_*(u)} \ge uN^d] \stackrel{(\ref{4.10})}{\ge} - c \,e^{-cu\,L^{d-2}} \,.
\end{align*}
As a result we find that:
\begin{equation}\label{4.52}
E[\Gamma_{c_{11}u}] - c\,e^{-cu \,L^{d-2}} \le E[\wt{\Gamma}_u] \le E[\Gamma_u] + c\,e^{-cu\,L^{d-2}}, \;\mbox{for $u > 0$}\,.
\end{equation}

\n
In the same fashion we also find that for $u > 0$,
\begin{equation}\label{4.53}
P[\Gamma_u < \wt{\Gamma}_u] \le c\,N^d \,e^{-c u\,L^{d-2}}, \;P[\wt{\Gamma}_{c^{-1}_{11} u} < \Gamma_u] \le c\,N^d \,e^{-cu\,L^{d-2}}\,.
\end{equation}

\n
Hence using the first inequalities in (\ref{4.52}) and (\ref{4.53}), we find that for $u,s > 0$:
\begin{equation*}
\begin{array}{l}
P[\Gamma_u < E[\Gamma_{c_{11} u}] - c \,e^{-c u\,L^{d-2}} -s ] \le P[\wt{\Gamma}_u < E[\wt{\Gamma}_u] - s] + c\,N^d \,e^{-c  u \,L^{d-2}}
\\[1ex]
\le \dis\frac{{\rm var}(\wt{\Gamma}_u)}{s^2}  + c \,N^d \,e^{- c u \,L^{d-2}}\,,
\end{array}
\end{equation*}

\n
and with (\ref{4.23}), optimizing over $r$, the claim (\ref{4.49}) follows. Using the rightmost inequalities of (\ref{4.53}) and of (\ref{4.52}), with $c^{-1}_{11} u$ in place of $u$, in the case of (\ref{4.52}), we analogously obtain (\ref{4.50}).
\end{proof}

\begin{remark}\label{rem4.4} \rm In the applications we discuss below, we will choose $L = [(\log N)^2]$, so that for given $u > 0$, and $N \ge c(u)$, 
\begin{equation}\label{4.54}
\sigma^2_{u,L,N} \le c\,u^{\frac{d}{d+1}} \;L^{\frac{d^2}{d+1}} \;N^{-\frac{d}{d+1}} \le c\,u^{\frac{d}{d+1}} \;(\log N)^{\frac{2d^2}{d+1}} \;N^{-\frac{d}{d+1}} \,,
\end{equation}
as follows from a straightforward upper bound of the expresson in (\ref{4.51}). \hfill $\square$
\end{remark}

We now turn to the first application of Theorem \ref{theo4.3} that sharpens (\ref{1.13}) into an estimate of the relative volume of the vacant set left by the walk at time $uN^d$.

\begin{corollary}\label{cor4.5} $(d \ge 3)$
\begin{equation}\label{4.55}
\lim\limits_N \;P[e^{-cu} \le \dis\frac{|E \backslash X_{[0,uN^d]}|}{N^d} \le e^{-c^\prime u}] = 1, \;\mbox{for  $u > 0$}\,.
\end{equation}
\end{corollary}

\begin{proof}
We choose $L = [(\log N)^2]$, and $\phi = \phi_0$, cf.~(\ref{4.2}), so that
\begin{equation}\label{4.56}
\Gamma_u \stackrel{(\ref{4.5})}{=} \dis\frac{1}{N^d} \;|E \backslash X_{[0,uN^d]}|, \;\mbox{for $u > 0$}\,,
\end{equation}
and with translation invariance
\begin{equation}\label{4.57}
E[\Gamma_u] = E[h (0,uN^d)] = P[0 \notin X_{[0,uN^d]}]\,.
\end{equation}

\n
Note that with the above choice for $L$, in view of (\ref{4.54}), $\sigma_{u,L,N}$ and $N^d \,e^{-cu\,L^{d-2}}$ tend to $0$ as $N$ tends to infinity. Choosing for instance $s = \sqrt{\sigma_{u,L,N}}$, the claim (\ref{4.55}) follows straightforwardly from (\ref{4.49}), (\ref{4.50}) and our estimates in (\ref{1.13}) on $E[\Gamma_{c_{11}u}]$ and $E[\Gamma_{c_{11}^{-1} u}]$.
\end{proof}

We recall that on the event $\cG_{\beta,t}$ defined in (\ref{2.50}), the vacant set left by the walk at time $t$ contains a well-defined unique giant component $O$, cf.~(\ref{2.51}), and $\cG_{\beta,uN^d}$ is typical under $P$ for large $N$, when $d \ge d_0$, and $u$ is small , cf.~(\ref{2.55}). As we will now see in this regime $O$ also typically occupies a non-degenerate fraction of the volume of $E$.

\begin{corollary}\label{cor4.6} $(d \ge d_0$, cf.~(\ref{2.39}))

\medskip
For any $\beta,\gamma \in (0,1)$, one has
\begin{equation}\label{4.58}
\lim\limits_N \;P\Big[\cG_{\beta,uN^d} \cap \Big\{\dis\frac{|O|}{N^d} \ge \gamma \Big\}\Big] = 1, \;\mbox{for small $u > 0$}\,.
\end{equation}
\end{corollary}

\begin{proof}
We choose $L = [(\log N)^2] \vee [c_0 \log N]$, cf.~Corollary \ref{cor2.5} and \ref{2.6}, and
\begin{align}
\phi(A) = 1\big\{& \mbox{for some $F \in \cL_2$, with $0 \in F$, 0 is connected to $S(0,L_0)$ in $F \backslash A\big\}$}, \label{4.59}
\\
& \mbox{for any $A \subseteq C(0) ( = B(0,L)$), with $L_0$ as in (\ref{2.51})}\,.   \nonumber
\end{align}

\n
In this case for large $N$ we have, cf.~(\ref{2.38}), (\ref{4.4}), (\ref{4.5}):
\begin{equation}\label{4.60}
\Gamma_u = \dis\frac{1}{N^d} \;\dsl_{x \in E} \;1_{\cC_{c_0,x,uN^d}} \stackrel{(\ref{2.54})}{\le} \dis\frac{|O|}{N^d}\,, \;\mbox{on} \;\cG_{\beta,uN^d}\;,
\end{equation}
and with translation invariance we find
\begin{equation}\label{4.61}
E[\Gamma_u] = P[\cC_{c_0,0,uN^d}]\,.
\end{equation}

\n
As already mentioned below (\ref{4.57}), $\sigma_{u,L,N}$ and $N^d \,e^{-c u \,L^{d-2}}$ tend to $0$ as $N$ tends to infinity. We can choose $s = \sqrt{\sigma_{u,L,N}}$ in (\ref{4.49}), so that 
\begin{equation*}
\lim\limits_N \;P\big[\Gamma_u \ge E[\Gamma_{c_{11}u}] - c e^{-cu L^{d-2}} - \sqrt{\sigma_{u,L,N}}\big] = 1\,.
\end{equation*}

\n
The claim (\ref{4.58}) then follows from (\ref{2.41}), (\ref{2.55}), and (\ref{4.60}).
\end{proof}

\begin{remark}\label{rem4.7} \rm ~

\medskip\n
1) When $d \ge d_0$, the above corollary shows that for small $u > 0$, when $N$ becomes large the giant component typically has non-degenerate volume in $E$. However it does not rule out the existence of other components in the vacant set with non-degenerate volume. Note that such components by the definition of the giant $O$, cf.~(\ref{2.51}), do not contain any connected sets $A \in \cA_2$ of $|\cdot |_\infty$-diameter $L_0 = [c_0 \log N]$ and in particular any segment of length $L_0$. 

\bigskip\n
2) When $d \ge 3$ and $u > 0$, the set visited by the walk up to time $uN^d$ typically constitutes a giant component as well. Indeed with Corollary \ref{cor4.5} it typically occupies a non-degenerate fraction of the volume of $E$, when $N$ is large. Moreover with a straightforward modification of (\ref{1.24}), see also Remark \ref{rem1.3}, we see that when $8  \wt{L} \le N$ and $u > 0$,
\begin{equation}\label{4.62}
\begin{split}
P[H_{B(0, \wt{L})} > uN^d] &\le c\exp \{- cu N^{d-2}\} + (1 - c( \wt{L}/N)^{d-2})_+^{k_* - 1}
\\
& \le c \exp\{- c u \wt{L}^{d-2}\}\,,
\end{split}
\end{equation}

\n
using the definition of $k_*$ below (\ref{1.12}). In particular choosing  $ \wt{L}=L_1 \stackrel{\rm def}{=} [c_{12} \,(\frac{\log N}{u})^{\frac{1}{d-2}}]$, we find that
\begin{equation}\label{4.63}
\lim\limits_N \,P[\mbox{for some $x$ in $E, X_{[0,uN^d]} \cap B(x,L_1) = \phi] = 0$, for all $u > 0$}\,.
\end{equation}

\n
So the set visited by the walk is ubiquitous as well, and typically comes within distance of order $(\log N)^{\frac{1}{d-2}}$ from any point of $E$.

\hfill $\square$

\end{remark}

\n
In fact holes in the vacant set of order $(\log N)^{\frac{1}{d-2}}$ do occur as well. More precisely, consider the maximal radius of an $|\cdot |_\infty$-ball contained in the vacant set at time $t$:

\begin{equation}\label{4.64}
\wh{L}(t) = \sup\{m \ge 0;   \mbox{for some $x$ in $E, X_{[0,t]} \cap B(x,m) = \phi$}  \}\,,
 \end{equation}
with the convention that $\wh{L}(t) =0$, when the right-hand side of (\ref{4.64}) is the empty set. 
 
\begin{corollary}\label{cor4.8} $(d \ge 3)$

\medskip
There exists $c_{13}$  $( < c_{12})$ such that
\begin{equation}\label{4.65}
\lim\limits_N \;P[L_2 \le \wh{L}(uN^d) \le L_1] = 1, \;\mbox{for  $u > 0$}\,,
\end{equation}
with $L_1$ defined above (\ref{4.63}) and $L_2 = [c_{13} \,(\frac{\log N}{u})^{\frac{1}{d-2}}]$.
\end{corollary}

\begin{proof}
In view of (\ref{4.63}) we only need to prove the lower bound. The argument uses a variation on the proof of Corollary \ref{cor4.5}. For $ \wt{L} \le  {(\log N)}^2 $ and large N, with a straightforward modification of (\ref{1.25}), we see that 
\begin{equation}\label{4.66}
\begin{split}
P[H_{B(0,\wt{L})} > uN^d] &\ge - c\exp \{- cu N^{d-2}\} + c(1 - c(\wt{L}/N)^{d-2})^{k^*}
\\
& \ge c \exp\{- c u \wt{L}^{d-2}\}\,,
\end{split}
\end{equation}
\n
using the definition of $k^*$ below (\ref{1.12}). Proceeding as in  Corollary \ref{cor4.5}, we then choose 
 $L = [(\log N)^2]$, and the monotone decreasing function $\phi$:
  \begin{equation*}
  \phi: A \subseteq C(0) \rightarrow \phi(A) = 1 \{B(0,\wt{L})\cap A = \phi\} \,.
 \end{equation*}
With this choice we find that 
 
\begin{equation*}
\Gamma_u \stackrel{(\ref{4.5})}{=} \dis\frac{1}{N^d} \;\dsl_{x \in E} \;1 \{B(x,\wt{L})\cap X_{[0,uN^d]} = \phi\} \;.
\end{equation*}

Setting $\wt{L} = [c_{13} \,(\frac{\log N}{u})^{\frac{1}{d-2}}]$  with $c_{13}$ small enough, we see with translation invariance that for large $N$:
\begin{equation}\label{4.67}
E[\Gamma_{c_{11}u}]  = P[H_{B(0,\wt{L})} >c_{11}uN^d]\ge N^{-\frac{1}{6}}\,.
\end{equation}
We then choose $s = \sqrt{\sigma_{u,L,N}}$ in (\ref{4.49}), and note with (\ref{4.54}) that for large $N$, $s = \sqrt{\sigma_{u,L,N}}$ is much smaller than $N^{-\frac{1}{6}}$, and that $\sigma_{u,L,N}$ and $N^d \,e^{-cu\,L^{d-2}}$ tend to $0$ as $N$ tends to infinity. As a result we obtain that
\begin{equation}\label{4.68}
\lim\limits_N \;P[\Gamma_u \le {\fr} N^{-\frac{1}{6}}] = 0\,,
\end{equation}
This is more than enough to prove the lower estimate in (\ref{4.65}). This concludes the proof of  Corollary\ref{cor4.8}.
\end{proof}

The above result exhibits a different asymptotic behavior from that of Bernoulli bond- (or site-) percolation on E, where for large $N$ the largest $|\cdot |_\infty$-ball contained in a cluster typically has size of order $(\log N)^{\frac{1}{d}}$, which is much smaller than $(\log N)^{\frac{1}{d-2}}$.


\begin{thebibliography}{10}

\bibitem{Aldo83}
D.~Aldous.
\newblock On the time taken by random walks on finite groups to visit every
  state.
\newblock {\em Z. Wahrscheinlichkeitstheorie verw. Gebiete}, 62(3):361--374,
  1983.

\bibitem{AldoFill99}
D.~Aldous and J.~Fill.
\newblock {\em Reversible Markov chains and random walks on graphs}.
\newblock http://www.stat.berkeley.edu/$\sim$aldous/RWG/book.html.

\bibitem{AlonBenjStac04}
N.~Alon, I.~Benjamini, and A.~Stacey.
\newblock Percolation on finite graphs and isoperimetric inequalities.
\newblock {\em Ann. Probab.}, 32(3A):1727--1745, 2004.

\bibitem{BrumHilh91}
M.J.A.M. Brummelhuis and M.J. Hilhorst.
\newblock Covering a finite lattice by a random walk.
\newblock {\em Physica A.}, 176(3):387--408, 1991.

\bibitem{DembPereRoseZeit04}
A.~Dembo, Y.~Peres, J.~Rosen, and O.~Zeitouni.
\newblock Cover times for {B}rownian motion and random walks in two dimensions.
\newblock {\em Ann. Math.}, 160(2):433--464, 2004.

\bibitem{DembPereRoseZeit06}
A.~Dembo, Y.~Peres, J.~Rosen, and O.~Zeitouni.
\newblock Late points for random walks in two dimensions.
\newblock {\em Ann. Probab.}, 34(1):219--263, 2006.

\bibitem{DembSzni06}
A.~Dembo and A.S. Sznitman.
\newblock On the disconnection of a discrete cylinder by a random walk.
\newblock {\em Probab. Theory Relat. Fields}, 136(2):321--340, 2006.

\bibitem{Durr91}
R.~Durrett.
\newblock {\em Probability: {T}heory and {E}xamples}.
\newblock Wadsworth and Brooks/Cole, Pacific Grove, 1991.

\bibitem{GrigTelc01}
A.~Grigoryan and A.~Telcs.
\newblock Sub-{G}aussian estimates of heat kernels on infinite graphs.
\newblock {\em Duke Math. J.}, 109(3):451--510, 2001.

\bibitem{HeydHofs05}
M.~Heydenreich and R.~van~der Hofstadt.
\newblock Random graph asymptotics on high dimensional tori.
\newblock {\em Preprint, {\rm arXiv:math.PR/0512522 v1}}, 2005.

\bibitem{Kest82}
H.~Kesten.
\newblock {\em Percolation theory for Mathematicians}.
\newblock Birkh\"auser, Basel, 1982.

\bibitem{Lawl91}
G.F. Lawler.
\newblock {\em Intersections of random walks}.
\newblock Birkh\"auser, Basel, 1991.

\bibitem{Lind92}
T.~Lindvall.
\newblock {\em Lectures on the coupling method}.
\newblock Dover Publications, Inc., New York, 1992.

\bibitem{Mont56}
E.W. Montroll.
\newblock Random walks in multidimensional spaces, especially on periodic
  lattices.
\newblock {\em J. Soc. Industr. Appl. Math.}, 4(4), 1956.

\bibitem{Szni06}
A.S. Sznitman.
\newblock How universal are asymptotics of disconnection times in discrete
  cylinders?
\newblock {\em Ann. Probab.},  36(1):1-53, 2008.

\bibitem{Szni98a}
A.S. Sznitman.
\newblock {\em Brownian motion, obstacles and random media}.
\newblock Springer, Berlin, 1998.

\end{thebibliography}
\end{document}